\documentclass[11pt]{article}

\usepackage{amsmath} 
\usepackage{amsthm} 
\usepackage{amssymb} 
\usepackage{enumitem}
\usepackage{hyperref}
\usepackage{color} % Voor het gebruik van kleuren

\newtheorem{thm}{Theorem}
\newtheorem{inspr}[thm]{}

\newenvironment{definitie}{\begin{itemize}\item[ ]\hspace{-26pt}\bf Definition \rm }{\end{itemize}}
\newenvironment{notatie}{\begin{itemize}\item[ ]\hspace{-26pt}\bf Notation \rm }{\end{itemize}}
\newenvironment{voorbeeld}{\begin{itemize}\item[ ]\hspace{-26pt}\bf Example \rm }{\end{itemize}}
\newenvironment{stelling}{\begin{itemize}\item[ ]\hspace{-26pt}\bf Theorem \rm }{\end{itemize}}
\newenvironment{propositie}{\begin{itemize}\item[ ]\hspace{-26pt}\bf Proposition \rm }{\end{itemize}}
\newenvironment{lemma}{\begin{itemize}\item[ ]\hspace{-26pt}\bf Lemma \rm }{\end{itemize}}
\newenvironment{opmerking}{\begin{itemize}\item[ ]\hspace{-26pt}\bf Remark \rm }{\end{itemize}}
\newenvironment{voorwaarde}{\begin{itemize}\item[ ]\hspace{-26pt}\bf Condition \rm }{\end{itemize}}
\newenvironment{probleem}{\begin{itemize}\item[ ]\hspace{-26pt}\bf Problem \rm }{\end{itemize}}

\renewcommand{\Bbb}{\mathbb} % hierdoor kan ik \Bbb blijven gebruiken

\newcommand{\defin}{\begin{inspr}\begin{definitie}}  %\def already defined
\newcommand{\edefin}{\end{definitie}\end{inspr}}

\newcommand{\notat}{\begin{inspr}\begin{notatie}}  %\not already defined
\newcommand{\enotat}{\end{notatie}\end{inspr}}

\newcommand{\voorb}{\begin{inspr}\begin{voorbeeld}}  %\not already defined
\newcommand{\evoorb}{\end{voorbeeld}\end{inspr}}

\newcommand{\stel}{\begin{inspr}\begin{stelling}}
\newcommand{\estel}{\end{stelling}\end{inspr}}

\newcommand{\prop}{\begin{inspr}\begin{propositie}}
\newcommand{\eprop}{\end{propositie}\end{inspr}}

\newcommand{\lem}{\begin{inspr}\begin{lemma}}
\newcommand{\elem}{\end{lemma}\end{inspr}}

\newcommand{\opm}{\begin{inspr}\begin{opmerking}}
\newcommand{\eopm}{\end{opmerking}\end{inspr}}

\newcommand{\voorw}{\begin{inspr}\begin{voorwaarde}}
\newcommand{\evoorw}{\end{voorwaarde}\end{inspr}}

\newcommand{\probl}{\begin{inspr}\begin{probleem}}
\newcommand{\eprobl}{\end{probleem}\end{inspr}}

\newcommand{\bew}{\vspace{-0.3cm}\begin{itemize}\item[ ] \bf Proof\rm: }
\newcommand{\bewo}{\vspace{-0.3cm}\begin{itemize}\item[ ] \bf Proof\rm (outline): }
\newcommand{\ebew}{\hfill $\qed$ \end{itemize}}

\newcommand{\ssnl}{\vskip 5pt} 
 
\newcommand{\nl}{\vskip 12pt}

\newcommand{\rood}{\color{red}}
\newcommand{\blauw}{\color{blue}}

\newcommand{\ot}{\otimes}

\newcommand{\tr}{\triangleright}

\newcommand{\tussenen}{\qquad\quad\text{and}\qquad\quad}

\newcommand{\inv}{^{-1}}

\numberwithin{thm}{section}  
\numberwithin{equation}{section} 

\newcommand{\keepcomment}[1]{}
\newcommand{\oldcomment}[1]{}

\newenvironment{niets}{\begin{itemize}[leftmargin=15pt] \item[ ]\hspace{- 15 pt}\bf   \rm }{\end{itemize}}
\newcommand{\nul}{\begin{inspr}\begin{niets}}
\newcommand{\enul}{\end{niets}\end{inspr}}

\begin{document}

\centerline{\bf \Large Algebraic quantum groups and duality II} % Verschillende andere mogelijkheden, zoals \large, \small, ...
\vspace{5pt}
\centerline{\bf \large Multiplier Hopf $^*$-algebras with positive integrals}
\vspace{13 pt}
\centerline{\it A.\ Van Daele \rm ($^*$)}
\vspace{20 pt}
{\bf Abstract} 
\nl
In the  paper \emph{Algebraic quantum groups and duality I} \cite{VD-part1}, we consider a pairing $(a,b)\mapsto\langle a,b\rangle$ of regular multiplier Hopf algebras  $A$ and $B$. When $A$ has integrals and when $B$ is the dual of $A$, we can describe the duality with an element $V$ in the multiplier algebra $M(B\ot A)$ satisfying and defined by $\langle V, a\ot b\rangle=\langle a,b\rangle$ for all $a,b$. Properties of the dual pair are formulated in terms of this multiplier $V$. It acts, in a natural way, on $A\ot A$ as the canonical map $T$, given by $T(a\ot a')=\Delta(a)(1\ot a')$.
\ssnl
In this second paper on the subject, we assume that the pairing is coming from  a  multiplier Hopf $^*$-algebra with positive integrals. In this case, the positive right integral on $A$ can be used to construct a Hilbert space $\mathcal H$. The duality $V$ now acts as a unitary operator on the Hilbert space tensor product $\mathcal H\ot\mathcal H$. This eventually makes it possible to complete the algebraic quantum group to a locally compact quantum group. 
\ssnl
The procedure to pass from the algebraic quantum group to the operator algebraic completion has been treated in the literature (see e.g.\ \cite{Ku-VD}) but the construction is rather involved because of the necessary use of left Hilbert algebras. 
\ssnl
In this paper, we  give a comprehensive, yet concise and somewhat simpler approach. It  should be considered as a \emph{springboard} to the more complicated theory of locally compact quantum groups.
\nl
Date: {\it 26 April 2023} 

\vskip 6cm
\hrule
\vskip 7 pt
\begin{itemize}
\item[($^*$)] Department of Mathematics, KU Leuven, Celestijnenlaan 200B,\newline
B-3001 Heverlee (Belgium). E-mail: \texttt{alfons.vandaele@kuleuven.be}
\end{itemize}

\newpage

\setcounter{section}{-1}  % Dit zorgt ervoor dat we met 0 beginnen voor de inleiding

\section{\hspace{-17pt}. Introduction} \label{sect:introduction}   % \input artikel0.tex%\newpage

In the introduction of the first paper on the subject \cite{VD-part1} we started by looking at a pair of {\it finite-dimensional} Hopf algebras $A$ and $B$ (over the field $\mathbb C$ of complex numbers), together with a non-degenerate pairing $(a,b) \mapsto \langle a,b \rangle$ from the Cartesian product $A\times B$ to $\mathbb C$. It is assumed that the product in $A$ induces the coproduct on $B$ and vice versa. Much of the information is encoded by means of the duality $V$, defined as an element in $B\ot A$ by $\langle V,a\ot b\rangle=\langle a,b\rangle$ where $a\in A$ and $b\in B$.
\ssnl 
We have some introductory sections in \cite{VD-part1} where we recall basic notions and properties of multiplier Hopf algebras, multiplier Hopf algebras with integrals and duality for those. We  discussed briefly more general pairings of regular multiplier Hopf algebras  $A,B$ (as in \cite{Dr-VD}), but mostly we worked with a pair  of multiplier Hopf algebras where $A$ is actually a regular multiplier Hopf algebra with integrals and where $B$ is the dual regular multiplier Hopf algebra with integrals. In this case, the duality $V$ turns out to belong to the multiplier algebra $M(B\ot A)$ of $B\ot A$. Again various results about $V$, proven in the finite-dimensional case, still holds here.
\ssnl
In this paper, we treat the case of a \emph{multiplier Hopf $^*$-algebra with positive integrals}. We will now use the term \emph{($^*$-)algebraic quantum groups} for such a multiplier Hopf $^*$-algebra. 
\nl
\bf Content of the paper\rm
\nl
In {\it Section} \ref{sect:prel} we recall some of the basic properties of multiplier Hopf $^*$-algebras with positive integrals. In particular, we review the analytic properties with its consequences. An important one is that the scaling constant $\tau$, defined by $\varphi\circ S^2=\tau \varphi$ where $\varphi$ is a left integral and $S$ the antipode, is in fact trivial. As a consequence a positive right integral exists if and only if there is a positive left integral. 
\ssnl
Further in this section, we discuss the Fourier transform, as studied already in the first paper,  see Section 2 of \cite{VD-part1}, and prove an analogue of Plancherel's formula. 
\nl
In {\it Section} \ref{sect:lcp} and further till {\it Section} \ref{sect:dual} we describe in detail how an algebraic quantum group $(A,\Delta)$ can be completed to a locally compact quantum group in the sense of Kustermans and Vaes (\cite{K-V1, K-V2, K-V3}). 
\ssnl
We start this procedure in Section \ref{sect:lcp} with the GNS-construction from the positive right integral. We obtain the duality $V$ on the Hilbert space level and we get the obvious von Neumann algebra $M$ containing $A$ as a dense $^*$-subalgebra with a coproduct $\Delta$ on $M$ that extends the original coproduct. 
\ssnl
The harder part in this procedure is to find the Haar weights, extending the integrals and proving that they are genuine invariant weights, making $(M,\Delta)$ into a locally compact quantum group. Some techniques and results from the theory of left Hilbert algebras are needed to do this. This is done in Section \ref{sect:hwt}. We begin with the completion of the right integral to the right Haar weight.
\ssnl
We have included an appendix with the most important notions and results about left Hilbert algebras and the associated construction of a weight, see {\it Appendix} \ref{sect:appA}. The more technical proofs are avoided in this appendix but are collected in yet another appendix, {\it Appendix} \ref{sect:appB}.
\ssnl
In {\it Section} \ref{sect:hwt}, we also need the square root of  the modular element $\delta$. This is used to construct the left Haar weight. Then the construction of the locally compact quantum group coming from an algebraic quantum group is complete. 
\ssnl
In the third paper in this series, \cite{VD-part3} we will look deeper into the analytic structure of an algebraic quantum group and we obtain various objects, like the modular automorphism $\sigma$ and $\sigma'$ in terms of Hilbert space operators and  some of the relations we found already in \cite{VD-part1} in a nice operator algebraic form. 
\ssnl
In {\it Section} \ref{sect:dual} we investigate the dual. We compare two possible approaches in this case. One possibility is to use the previous procedure for the dual $\widehat A$ of $A$ in the sense of algebraic quantum groups. The other one is to use the construction of the dual of the locally compact quantum group in the sense of locally compact quantum groups. We show that these give the same dual locally compact quantum group as expected.
\ssnl
In the last section, {\it Section} \ref{sect:concl}, we draw some conclusions and discuss possible further research.
\nl
\bf Motivation \rm
\nl
The material that we treat in this note has been studied in earlier papers, first in \cite{Ku-VD} and later, partly also in \cite{DC-VD}. However, the approach here, certainly when compared with the first paper \cite{Ku-VD}, is simpler and more direct. We need some results of left Hilbert algebras and the relation with normal semi-finite weights, but we try to take advantage of the concrete situation where we have to use them. Adding an appendix about the construction of a normal semi-finite weight associated with a left Hilbert algebra should make this easier for the reader.
\nl
\bf Notations and conventions \rm
\nl
We only work with $^*$-algebras over $\Bbb C$. A linear functional $\omega$ on a $^*$-algebra $A$ is called positive if $\omega(a^*a)\geq 0$ for all $a\in A$. It is  faithful if $\omega(a^*a)=0$ is only true when $a=0$.
\ssnl
We use the tensor product symbol in various cases. When we have a  $^*$-algebra $A$ (without any topology), then $A\ot A$ denotes the algebraic tensor product. It is again a $^*$-algebra in a natural way. When we have a $C^*$-algebra $P$ we consider the minimal $C^*$-tensor product $P\ot P$ and when we have a von Neumann algebra $M$, we consider the von Neumann algebra tensor product $M\ot M$. Finally, $\mathcal H\ot \mathcal H$ denotes the Hilbert space tensor product if $\mathcal H$ is a Hilbert space. It should always be clear from the context which type we are using. If not, it will be explicitly stated.
\ssnl
For a coproduct $\Delta$ on a $^*$-algebra $A$ (without topology), as we use it in this theory, we assume that $\Delta(a)(1\ot a')$ and $(a\ot 1)\Delta(a')$ are in $A\ot A$ for all $a,a'\in A$. For a coproduct on an operator algebra, a C$^*$-algebra or a von Neumann algebra, we do not have such a restriction. In Section \ref{sect:lcp}, we recall the notion of a coproduct in these two cases.
\ssnl

We sometimes use the {\it Sweedler notation} for a coproduct. See the item \emph{Notations and conventions} in the Introduction of Part I \cite{VD-part1} for more information about this.
\ssnl
Finally, in order to avoid too many different notations, subscripts, etc.,  we use the same symbol for different objects. 
We believe this will not lead to confusions. In any case, whenever there is some possible doubt, we will be more explicit.
\nl
\bf More about conventions \rm
\nl
It is well known that different conventions are used in this theory. Roughly speaking, there are two basic types. There is a difference between the purely algebraic treatment and the operator algebraic versions. But also within the operator algebraic approach, different conventions are used. This is very unfortunate as it sometimes makes it difficult to compare the various formulas. They then all look slightly different. Also, there is no systematic {\it dictionary} to translate from one to the other.
\ssnl
In the purely algebraic theories, the coproduct on the dual algebra $B$ satisfies\\
$\langle a\ot a',\Delta(b)\rangle=\langle aa',b\rangle$ when $a,a'\in A$  and $b\in B$. 
On the other hand, in the operator algebraic approaches (say the theory of locally compact quantum groups), it is common to take the flipped coproduct on the dual, that is to assume that $\langle a\ot a',\Delta(b)\rangle=\langle a'a,b\rangle$ when $a,a'\in A$ and $b\in B$, (whenever such a formula makes sense).
\ssnl
In this note, we will always use the {\it algebraic convention}, also when we eventually come to the construction of a locally compact quantum group in Section \ref{sect:lcp}.
\ssnl
Another difference is mostly encountered only within the operator algebraic approach. Essentially in one case the duality is seen as a unitary operator as we construct it with the right integral, i.e.\ the right regular representation, while in the other case, the duality is the unitary operator we get from the left integral (by applying the Fourier transform and then interchanging the roles of $A$ and the dual $\widehat A$, that is the left regular representation).
\ssnl
In this paper, we will work with the {\it right regular representation} and we will try to explain why this is a better and more natural choice. We refer to earlier discussions on this matter in Section \ref{sect:concl} as well. \oldcomment{NIet vergeten. \rood Check!}{}
\ssnl
Finally, again within the theory of locally compact quantum groups, different authors use different conventions about the polar decomposition of the antipode. We will explain this in \cite{VD-part3} when this decomposition is introduced. We will then mention which of the two choices and why we do so. \oldcomment{Niet vergeten. \rood Check!}{}
\ssnl
We have a similar discussion about these different conventions in \cite{VD-warsaw}.
\nl
\bf Basic references \rm
\nl
For the theory of Hopf algebras, we refer to the well-known books by Abe \cite{Ab} and Sweedler \cite{Sw}. See also the more recent work by Radford \cite{R-bk}. The original work on multiplier Hopf algebras is \cite{VD-mha} and for multiplier Hopf algebras with integrals, it is \cite{VD-alg}. The use of the Sweedler notation for multiplier Hopf algebras has been explained in e.g.\ \cite{VD-tools} and more recently in \cite{VD-sw}.
\ssnl
Pairings of multiplier Hopf algebras have been first studied in \cite {Dr-VD}. Actions of multiplier Hopf algebras are studied in \cite{Dr-VD-Z}.
\ssnl
Locally compact quantum groups have been considered by various authors and there are also different approaches. Our main references here are \cite{K-V2} and \cite{K-V3}. See also \cite{VD-sigma} for a slightly simplified treatment of the theory. The reader is also advised to look at the more recent survey paper \cite{VD-warsaw}.
\ssnl
For the locally compact quantum groups arising from algebraic quantum groups, the first paper to consider is \cite{Ku-VD}. Related is the analytic structure of an algebraic quantum group, see \cite{Ku}.  A more recent treatment is found in \cite{DC-VD}. %Finally, in 
\nl%\nl
\bf Acknowledgments \rm
\nl
I am very grateful to M.B.\ Landstad and other colleagues and friends, both at the University of Trondheim and the University of Oslo (where part of these notes were developed) for the nice and fruitful atmosphere during my regular visits to these departments.

\section{\hspace{-17pt}. Multiplier Hopf $^*$-algebras with positive integrals} \label{sect:prel}% \input artikel1.tex %\newpage

n the previous paper \cite{VD-part1} we have collected some formulas and results about regular multiplier Hopf algebras wiht integrals. In the first section of \emph{this paper}, we add results for the case where these multiplier Hopf algebras are \emph{multiplier Hopf $^*$-algebras}. 
Mostly, we will work with integrals that are positive but we start without this requirement.
\ssnl 
Here is an example with integrals, but no positive integrals.

\voorb
 Let $\lambda$ be a complex number  that is a root of unity, different from $1$ and $-1$. Let $n$ be the smallest natural number so that $\lambda^{2n}=1$. Consider the unital $^*$-algebra generated by {\it self-adjoints} element $a,b$ where 
\begin{equation*}
a \text{ is invertible}, \qquad b^n=0, \qquad  ab=\lambda ba.
\end{equation*}
 Remark that these three conditions do not conflict with the self-adjointness of the elements $a$ and $b$. This is true for the third condition because $\lambda$ has modulus $1$. 
  \ssnl
This algebra is a Hopf $^*$-algebra with a coproduct determined by 
\begin{equation*}
\Delta(a)=a\ot a
\tussenen
\Delta(b)=a\ot b + b\ot a^{-1}.
\end{equation*}
See Proposition 5.6 in \cite{VD-alg}. 
It is shown in Proposition 5.8 of \cite{VD-alg} that this Hopf $^*$- algebra has integrals. They can not be chosen to be positive as these integrals are $0$ in $1$.
\evoorb

Remark that $b$ is self-adjoint and $b^n=0$. This can not happen in an operator algebra.
\ssnl
We could add the condition that $a^{2n}=1$ as this is compatible with $ab=\lambda ba$ because $\lambda^{2n}=1$. Then we obtain such a finite-dimensional example. See Section 3 in \cite{VD-fd}. \oldcomment{Reference?}{}
\nl
Next we claim that there is always a self-adjoint integral if there is an integral. We use $\overline\omega$ for the functional on $A$ given by $\overline\omega(a)=\overline{\omega(a)}$ whenever $\omega$ is a linear functional on a $^*$-algebra $A$.

\prop\label{prop:1.2a}
Let $\varphi$ be a left integral, then $\overline\varphi$ is also a left integral and $\varphi+\overline\varphi$ is a self-adjoint left integral.
\eprop
\bew
For all $a\in A$ we have $\Delta(a^*)=\Delta(a)^*$ because $\Delta$ is assumed to be a $^*$-map. If we apply $\varphi$ on the second leg we obtain
\begin{equation*}
\varphi(a^*)1=(\iota\ot\varphi)(\Delta(a)^*)=((\iota\ot\overline\varphi)(a))^*
\end{equation*}
and so $(\iota\ot\overline\varphi)(a)=\overline\varphi(a)1$ because $1^*=1$.
\ebew

Remark that, by the uniqueness of left integrals, we will have that both $\overline\varphi$ and $\varphi+\overline\varphi$ are scalar multiples of the original left integral $\varphi$.
\nl
Next suppose that $\varphi$ is a positive left integral.  Now we have the following result. Recall that the scaling constant is defined by $\varphi\circ S^2=\tau\varphi$.

\prop
The right integral $\psi$, defined as $\varphi\circ S$, is self-adjoint if and only if the scaling constant  $\tau$ is equal to $1$.
\eprop

\bew For all $a\in A$ we have
\begin{align*}
\overline\psi(a)&=\psi(a^*)^-=\varphi(S(a^*))^-\\
&=\varphi(S^{-1}(a)^*)^-=\varphi(S^{-1}(a))\\
&=\tau^{-1}\varphi(S(a))=\tau^{-1}\psi(a).
\end{align*}
We use $\overline\lambda$ for the complex conjugate of a complex number $\lambda$ and $\lambda^-$ as a short hand notation for $\overline\lambda$.
Therefore $\overline\psi=\tau\inv\psi$. This implies the result
\ebew

It can be shown that for an algebraic quantum group, when $\varphi$ is a positive left integral, then $\varphi\circ S$ is again positive. This result is not obvious and it was an open problem for some time. It was first proven in \cite{DC-VD}. As a consequence of the previous proposition, the scaling constant has to be $1$.  
\ssnl
We will give a short proof of these important results. We begin however with adding some formulas to the collection obtained in \cite{VD-part1}  in the case of a multiplier Hopf $^*$-algebra with integrals (not necessarily positive). 
\nl
\bf Some basic formulas in the involutive case \rm
\nl
From the uniqueness of the antipode, it follows that $a\mapsto S(a)^*$ is involutive, i.e.\ that $S(S(a)^*)^*=a$ for all $a$. We also have that $S(aa')^*=S(a)^*S(a')^*$ for all $a,a'$.
\ssnl
Next are some properties of the other objects in the case of a multiplier Hopf $^*$-algebra with integrals, not necessarily positive, that are relatively easy to obtain.

\prop 
Let $(A,\Delta)$ be a multiplier Hopf $^*$-algebra with integrals.
Then we get the following properties of the associated objects.
\begin{itemize}[noitemsep] 
\item[i)]
The scaling constant $\tau$ must have modulus $1$.
\item[ii)]
The modular element $\delta$ is self-adjoint.
\item[iii)]
For the modular automorphisms $\sigma$ and $\sigma'$ of respectively the left integral and the right integral, we have for all $a\in A$ that 
\begin{equation}
\sigma(a^*)=\sigma^{-1}(a)^* 
\tussenen 
{\sigma'}(a^*)={\sigma'}^{-1}(a)^*. \label{eqn:1.1}
\end{equation}
\end{itemize}
\eprop 

\bew
To prove these results, we need to work with self-adjoint integrals. As we have seen in Proposition \ref{prop:1.2a}, they always exist when there is an integral.   
\ssnl
i) Now take a self-adjoint left integral $\varphi$. By the definition of the scaling constant we have
$\varphi(S^2(a))=\tau \varphi(a)$ (see Proposition A.3 in \cite{VD-part1}).
Take complex conjugates, use that $\varphi$ is self-adjoint and that $S^2(a)^*=S^{-2}(a^*)$. Then we find
$\varphi(S^{-2}(a^*))=\overline\tau\varphi(a^*)$ for all $a$. It follows that $\tau^{-1}=\overline\tau$ and hence $|\tau|=1$.
\vskip 3pt
ii) By taking adoints of the two sides of the equality $(\varphi\ot\iota)\Delta(a)=\varphi(a)\delta$, obtained in Proposition  1.6 in \cite{VD-part1}, we find $\delta^*=\delta$. \oldcomment{\rood Aanvullen}{}
\vskip 3pt
iii) The formulas (\ref{eqn:1.1}) are also immediate consequences of the defining formulas for the modular automorphisms in Proposition 1.7 in \cite{VD-part1} if we take self-adjoint integrals. 

\ebew

Remark that the Equalities (\ref{eqn:1.1}) in item iii) will imply e.g.\  that the two equalities in Proposition 1.9 of \cite{VD-part1} will be consequences of each other. In fact, one can verify that the formulas we have in Section 1 of \cite{VD-part1} are consistent with the ones of the proposition above. \oldcomment{Relevantie?}{}
%\ssnl 
\oldcomment{\rood Check this item once more}
\nl
\bf More properties when there is a positive left integral \rm
\nl
From now on, we assume that $(A,\Delta)$ is a multiplier Hopf $^*$-algebra with a {\it positive left integral $\varphi$}.
\ssnl
One can show that $\delta$ must be a positive element in $M(A)$. This is explained in the following proposition. Here, an element in $A$ is called positive if it is of the form $a^*a$ for some $a$. 

\prop
The element $a^*\delta a$  is a positive element in $A$ for all $a$.
\eprop

\bew
Take any  $b$ in $A$ and write $\Delta(b)(1\ot a)=\sum_i q_i\ot p_i$. Then
\begin{equation*}
(1\ot a^*)\Delta(b^*b)(1\ot a)=\sum_{i,j} q_i^*q_j \ot p_i^*p_j
\end{equation*}
and if we apply $\varphi$ on the first leg, we get (by using Proposition 1.6 of \cite{VD-part1}
\begin{equation*}
\varphi(b^*b)a^*\delta a=\sum_{i,j} \varphi(q_i^*q_j)p_i^*p_j 
\end{equation*}
As $\varphi$ is positive, we have a positive matrix with entries $\varphi(q_i^*q_j)$. We can diagonalize it and this will eventually lead to an equality of the form $\varphi(b^*b)a^*\delta a=c^*c$ where $c$ is a linear combination of the elements $(p_i)$.
\ebew

Although that is expected, it is not clear wether also $\delta^{-1}$ is positive. But indeed it will follow later. We could use the above argument with the right integral, but then we would need to know already that also the right integral is positive.

\ssnl
We now proceed to argue that this is indeed the case.
\ssnl
First we show the following. The proof is inspired by the arguments given in \cite{DC-VD} (see e.g.\ Proposition 2.8 in \cite{DC-VD}).

\prop \label{prop:1.6}
Assume that $(A,\Delta)$ is a multiplier Hopf $^*$-algebra with a positive left integral $\varphi$. Then for all elements $a\in A$, there is a finite-dimensional subspace containing the elements $\delta^n a$ and $a\delta^n$ where $n\in \mathbb Z$.
\eprop

\bew
i) Let $a,b\in A$ and assume $b\neq 0$. We write $a\ot b^*b$ as a finite sum
\begin{equation*}
a\ot b^*b=\sum_i \Delta(p_i)(q_i\ot 1)
\end{equation*}
where $(p_i)$ and $(q_i)$ are elements in $A$. Multiply from the left with $\delta^n\ot\delta^n$ and use that $\Delta(\delta)=\delta\ot\delta$ (see Proposition 1.6 in \cite{VD-part1}). We find that 
\begin{equation*}
\delta^n a\ot \delta^n b^*b=\sum_i \Delta(\delta^np_i)(q_i\ot 1)
\end{equation*}
for all $n\in\mathbb Z$. Now apply the left integral $\varphi$ on the second leg of this equation. We see that 
\begin{equation}
\varphi(\delta^n b^*b)\delta^n a=\sum_i\varphi(\delta^np_i)q_i.\label{eqn:C.2}
\end{equation}
If $n=2m$ for $m\in \mathbb Z$ we have
\begin{equation*}
\varphi(\delta^n b^*b)=\varphi(\delta^{2m}b^*b)=\varphi(\delta^m b^*b\sigma(\delta^m))=\tau^{-m}\varphi(\delta^m b^*b\delta^m).
\end{equation*}
We have used that $\sigma(\delta)=\tau^{-1}\delta$, see Proposition 1.9 in \cite{VD-part1}.\oldcomment{\rood Aanvullen.}{}
\vskip 3pt
Because $\delta^*=\delta$, it follows from the faithfulness of $\varphi$ that this can only be $0$ if $b\delta^m=0$. However, by assumption $b\neq 0$ so that $\varphi(\delta^n b^*b)\neq 0$, at least when $n$ is even. Then we see from Equation \ref{eqn:C.2} that $\delta^n a$ belongs to the space spanned by the elements $(q_i)$, at least for $n$ even. 
\ssnl
ii) If we apply this for $a$ as well as for $\delta a$, we get that $\delta^na$ belongs to a finite-dimensional subspace of $A$. In a similar way, or by taking adjoints we get a finite-dimensional space containing all the elements  $a\delta^n$.
\ebew

The following is now an easy consequence.

\prop\label{prop:1.7a}
Let $(A,\Delta)$ be a multiplier Hopf $^*$-algebra with positive integrals. Then $A$ is spanned by elements that are simultaneously eigenvectors for left and right muliplication by $\delta$. Moreover, all eigenvalues are positive.
\eprop
\bew
i) Fix $a$ in $A$ and let $L$ be the space spanned by the elements $\delta^n a$ with $n\in \mathbb Z$. The space is left invariant by multipication from the left with $\delta$. There is a scalar product on $L$ given by $\langle x', x \rangle =\varphi(x^*x')$. And because $\delta$ is self-adjoint, it acts as a self-adjoint operator on this finite-dimensional inner product space. It can be diagonalized, with real eigenvalues. In fact, because $a^*\delta a$ is positive, these eigenvalues have to be positive.
\ssnl
ii) We have the same property for right multiplication with $\delta$. The result can be obtained in a similar way, or by taking adjoints.
\ssnl
iii) Because left and right multiplication by $\delta$ commute, we can get vectors that are simultaneously eigenvectors for left and right multiplication.
\ebew

We use this result now. See Theorem 3.4 in \cite{DC-VD}.

\prop\label{prop:1.7}
If $(A,\Delta)$ is a multiplier Hopf $^*$-algebra with a positive integral, then the scaling constant $\tau$ is equal to $1$.
\eprop 

\bew
Take any non-zero element $a$ with the property that $\delta a=\lambda a$ for some $\lambda\in \mathbb R$. We have seen that such elements exist as a consequence of the previous lemma. Then 
\begin{align*}
\lambda \varphi(aa^*)
&=\varphi(\delta aa^*)
=\varphi(aa^*\sigma(\delta))\\
&=\tau\inv\varphi(aa^*\delta)
=\tau\inv\varphi(a(\delta a)^*)\\
&=\tau\inv\lambda\varphi(aa^*).
\end{align*}
This implies that $\tau=1$.
\ebew

One could also argue as follows. We have $\sigma(\delta a)=\tau^{-1}\delta\sigma(a)$ because $\sigma(\delta)=\tau\inv\delta$ (see Proposition 1.9 in \cite{VD-part1}). Since all eigenvalues of left multiplication with $\delta$ have to be posisitve, we will need $\tau$ positive. But we know already that it has modulus $1$, hence $\tau=1$. 
\vskip 3pt
We now show that $\varphi\circ S$ is positive when $\varphi$ is positive. See Corollary 3.6 in \cite{DC-VD}.

\stel\label{stel:1.8}
Let $(A,\Delta)$ be multiplier Hopf $^*$-algebra with a positive left integral $\varphi$. Then $\varphi\circ S$ is a positive right integral.
\estel

\bew
Take an element $a$ and write it as a sum $\sum_i a_i$ where $a_i$ are elements satisfying $ a_i\delta=\lambda_i a_i$ for some $\lambda_i\in \mathbb R$. This is possible because Proposition \ref{prop:1.6}. We assume that all these eigenvalues are different from each other. For all  $j,k$ we find
\begin{equation*}
\varphi( a_j^* a_k\delta)=\lambda_k\varphi(a_j^*a_k)
\tussenen
\varphi(\delta a_j^* a_k)=\lambda_j\varphi(a_j^*a_k).
\end{equation*}
These two expressions are the same as $\sigma(\delta)=\delta$.
If $j\neq k$ we have $\lambda_j\neq \lambda_k$ and it follows that $\varphi(a_j^*a_k)=0$.
\vskip 3pt
Then 
\begin{equation*}
\varphi(a^*a\delta)=\sum_{j,k}\lambda_k\varphi(a_j^*a_k)=\sum_k\lambda_k\varphi(a_k^*a_k).
\end{equation*}
This is positive as $\lambda_k$ and $\varphi(a_k^*a_k)$ are both positive for all $k$. 
\ssnl
As $\varphi(S(a))=\varphi(a\delta)$, the result follows.
\ebew

As we have seen, this also implies that the scaling constant is $1$. However, we can not use this argument because in the proof, we have used that the scaling constant is $1$.
\keepcomment{
Ook een voorbeeld geven met een niet unimodulaire totaal onsamenhangende groep om een en ander met die $\delta$ te illustreren? \rood Wat heb ik hier mee bedoeld?}{}
\nl
\bf The analytic structure of an algebraic quantum group\rm
\nl
It was shown by Kustermans in \cite{Ku} that the maps $S^2$, $\sigma$, $\sigma'$ as well left and right multiplication with $\delta$ have {\it analytical extensions}. 
\ssnl
An argument like the one given in the proof of Proposition  \ref{prop:1.6}, can be used to obtain such results in a simpler way. See Section 3 in \cite{DC-VD}. In particular, there is Theorem 3.5 of \cite{DC-VD}:

\stel
Let $(A,\Delta)$ be a multiplier Hopf $^*$-algebra with positive integrals. Then $A$ is spanned by elements which are simultaneously eigenvectors for $S^2$, $\sigma$, $\sigma'$ and left and right multiplication by $\delta$. Moreover, all the eigenvalues are strictly positive.
\estel

We have shown already part of this result in Proposition \ref{prop:1.7}. For a proof of the other part, we refer to Section 2 of \cite{VD-part3} where we say more about the analytic sturcture of a $^*$-algebraic quantum group.
We will not need this here. 
\ssnl
We will use the result of Proposition \ref{prop:1.7}  here to give an other argument for the positivity of the right integral $\varphi\circ S$ that we obtained in Theorem \ref{stel:1.8}. This is also found in Section 3 of \cite{DC-VD}. The following result  will also be needed in Section \ref{sect:dual} when we construct the left Haar weight.

\prop\label{prop:1.11a}
There is a self-adjoint multiplier $\gamma$ in $M(A)$ satisfying $\gamma^2=\delta$. It will be the unique element satisfying this property if we require it to be positive. In that case, we have $\sigma(\gamma)=\gamma$.
\eprop

\bew We define a multiplier $\gamma$ by
\begin{align*}
\gamma a&=\lambda^{\frac12} a \text{ if } a\in A \text{ and }\delta a=\lambda a \\
a\gamma &=a\lambda^{\frac12} \text{ if } a\in A \text{ and } a\delta =\lambda a.
\end{align*}
One can verify that this is well-defined in $M(A)$. To do this, we use that any element in $A$ is a sum of eigenvectors. 
\ssnl
Now assume  that  $\sum_i a_i=0$ where for each $i$ we have $\delta a_i=\lambda_i a_i$ for some real number $\lambda_i$.  Assume that all these eigenvalues are different from each other. For the scalar product induced by $\varphi$, we have that the elements $a_i$ are mutually orthogonal. Hence $a_i=0$ for all $i$. Then also $\sum_i \lambda^{\frac12}a_i=0$. This means that we can define the map $a\mapsto \gamma a$. Similarly for the other map.
\ssnl
We have a multiplier of $A$ because $a'(\gamma a)=(a'\gamma)a$ for all eigenvectors and hence for all elements $a,a'$.
\ssnl
It is clear that $\gamma^2=\delta$ and that it is the unique element satisfying this property having positive eigenvalues.
\ssnl
Finally, from the uniqueness and because $\sigma(\delta)=\delta$, it also follows that $\sigma(\gamma)=\gamma$.

\ebew

We denote this element $\gamma$ by $\delta^{\frac12}$.
\ssnl
In fact, we can define $\delta^{z}$ for any complex number $z$ like that. And we get $\delta^z\delta^{z'}=\delta^{z+z'}$ for all $z,z'\in \mathbb C$. Moreover, when $z=1$ we get the original multiplier $\delta$ while for $z=\frac12$ we get the multiplier $\gamma$ of the proposition. 
\ssnl
This result can be used to give another proof of Theorem \ref{stel:1.8}. Indeed, we have
\begin{equation*}
\varphi(S(a))=\varphi(a\delta)=\varphi(\delta^{\frac12}a\delta^{\frac12})
\end{equation*}
for all $a$. Because $\varphi$ is assumed to be positive and $\delta^{\frac12}$ is self-adjoint, we get that also $\varphi\circ S$ is positive. This is the argument given in the proof of Corollary 3.6 of \cite{DC-VD}.

%\newpage %%%%%%%%%
\nl
\bf  The dual and the Fourier transform \rm
\nl
We have considered the dual of a multiplier Hopf algebra with integrals in Section 2 of \cite{VD-part1}. In the case of a multiplier Hopf $^*$-algebra with integrals, the dual is again a multiplier Hopf $^*$-algebra if we define $b^*$ for $b\in \widehat A$ by 
\begin{equation*}
\langle a,b^*\rangle=\langle S(a)^*,b\rangle^-.
\end{equation*}
If the left integral $\varphi$ on $A$ is positive, then the right integral $\widehat\psi$, defined on $\widehat A$ as in Proposition 2.1 of \cite{VD-part1}, is again positive. This result is already found in Proposition 4.8 of \cite{VD-alg}. We give a proof of this result below.
\nl
In Section 2 of \cite{VD-part1}, we have two forms of the Fourier transform from $A$ to $\widehat A$. 
% \ssnl
On the one hand there is the map $a\mapsto \varphi(\,\cdot\,a)$
 as in Proposition 2.6 of \cite{VD-part1}. 
On the other hand we have $a\mapsto \psi(S(\,\cdot\,)a)$  as in Proposition 2.7 of \cite{VD-part1}. 
\ssnl
They behave nicely with respect to the involutive structure as we see in the next proposition.

\prop \label{prop:C.10}
i) Define $\widehat \psi(b)=\varepsilon(a)$ if $b=\varphi(\,\cdot\, a)$.  Then $\widehat \psi(b^*b)=\varphi(a^*a)$ when again  
$b=\varphi(\,\cdot\, a)$.
\vskip 3pt
ii) Define $\widehat\varphi(b)=\varepsilon(a)$ when $b=\psi(S(\,\cdot\,)a)$ for $a\in A$. Then $\widehat\varphi(b^*b)=\psi(a^*a)$ when again $b=\psi(S(\,\cdot\,)a)$ for $a\in A$.
\eprop

\bew
i) Take $a\in A$ and define $b=\varphi(\,\cdot\, a)$. Take $c\in A$. Then
\begin{align*}
\langle c, b^*b \rangle 
&=\sum_{(c)}\langle c_{(1)}, b^*\rangle \langle c_{(2)},b \rangle \\
&=\sum_{(c)}\langle c_{(1)}, b^*\rangle\varphi (c_{(2)}a) \\
&=\sum_{(a)}\langle S^{-1}(a_{(1)}), b^*\rangle\varphi (ca_{(2)}).
\end{align*}
We have used that $(\iota\ot\varphi)(1\ot c)\Delta(a))=S((\iota\ot\varphi)(\Delta(c)(1\ot a)))$ (see Equation (1.2) in Proposition 1.1 of \cite{VD-part1}. By the definition of $\widehat\psi$ we get
\begin{equation*}
\widehat\psi(b^*b)=\sum_{(a)}\langle S^{-1}(a_{(1)}), b^*\rangle\varepsilon(a_{(2)})=\langle S^{-1}(a), b^*\rangle=\langle a^*,b\rangle^-.
\end{equation*}
This proves that $\widehat \psi(b^*b)=\varphi(a^*a)$.
%\vskip 3pt
\ssnl
ii) Take $a\in A$ and define $b=\psi(S(\,\cdot\,)a)$. Take $c\in A$. Then
\begin{align*}
\langle c,b^*b \rangle
&=\sum_{(c)}\langle c_{(1)}, b^*\rangle \langle c_{(2)},b \rangle \\
&=\sum_{(c)}\langle c_{(1)}, b^*\rangle \psi(S(c_{(2)})a) \\
&=\sum_{(Sc)}\langle S^{-1}((Sc)_{(2)}), b^*\rangle \psi((Sc)_{(1)}a) \\
&=\sum_{(a)}\langle S^{-1}(S(a_{(2)}), b^*\rangle \psi(S(c)a_{(1)})\\
&=\sum_{(a)}\langle (a_{(2)}, b^*\rangle \psi(S(c)a_{(1)}).
\end{align*}
We have used Equation 1.3 of Proposition 1.1 of \cite{VD-part1}
\begin{equation*}
S((\psi\ot\iota)((c'\ot 1)\Delta(a)))=(\psi\ot\iota)(\Delta(c')(a\ot 1))
\end{equation*}
for $c'=S(c)$. The we get
\begin{equation*}
b^*b=\sum_{(a)}\langle (a_{(2)}, b^*\rangle \psi(S(\,\cdot\,)a_{(1)}).
\end{equation*}
and 
\begin{align*}
\widehat\varphi(b^*b)
&=\sum_{(a)}\langle (a_{(2)}, b^*\rangle \varepsilon(a_{(1)})\\
&=\langle a,b^*\rangle=\langle S(a)^*,b\rangle^- \\
&=\psi(S(S(a)^*)a)^-=\psi(a^*a)^-.
\end{align*}
This proves that also here we have $\widehat\varphi(b^*b)=\psi(a^*a)$.
\ebew

We see from these results that a positive right integral on the dual $B$ exists if a positive left integral on $A$ exists. Also a positive left integral exists on $B$ if a positive right integral exists on $A$. By the previous results we get that we have positive left and right integrals on $A$ and on $B$ as soon as we have a positive left integral on $A$.

\opm
In fact, there was no need to prove the second part of the previous proposition. 
It follows by duality because the inverse of the map $a\mapsto \varphi(\,\cdot\,a))$ is the map $b\mapsto \widehat\psi(\widehat S(\,\cdot\,)a)$ (see Proposition 2.6 of \cite{VD-part1}) 
while also the inverse of $a\mapsto \psi(S(\,\cdot\,)a)$ is $b\mapsto \widehat\varphi(\,\cdot\, b)$ (see Proposition 2.7 in \cite{VD-part1}).%\mycomment{Find concrete references to these results in \cite{VD-part1}. }{}
\eopm

\section{\hspace{-17pt}. The associated locally compact quantum group}\label{sect:lcp} % \input artikel2a.tex %\newpage

In \cite{VD-part1} we have started with a dual pair of multiplier Hopf algebras. We have considered the Heisenberg commutation relations induced by the pairing and  studied the Heisenberg algebra. We quickly moved to the more restrictive case arising from a regular multiplier Hopf algebra $A$ with integrals. The dual $\widehat A$ was taken for $B$ and the duality was the natural duality between $A$ and its dual $\widehat A$. In that case, the duality is an invertible element $V$ in the multiplier algebra $M(B\ot A)$ with nice properties reflecting the duality.
\ssnl
Here we assume that the multiplier Hopf algebra $A$ is a {\it multiplier Hopf $^*$-algebra}, that it has \emph{integrals} and that the integrals are {\it positive}. We know that the dual $\widehat A$ is again a multiplier Hopf $^*$-algebra with positive integrals. The natural duality between the two provides a dual pair of multiplier Hopf $^*$-algebras. Indeed, we have the following result.

\prop
The pairing between a multiplier Hopf $^*$-algebra $A$ with positive integrals and its dual $\widehat A$ is a pairing of multiplier Hopf $^*$-algebras in the sense of Definition 2.8 in \cite{Dr-VD}.
\eprop

The result is mentioned in \cite{Dr-VD} but without a detailed proof. In Proposition 3.7 of \cite{VD-part1} we have formulated and proved this result for regular multiplier Hopf algebras with integrals. We refer to the proof of that result in \cite{VD-part1}. Then the only thing that is missing to complete the proof of the above proposition is the basic relation defining the involultion on the dual. It says that $\langle a,b^* \rangle=\langle S(a)^*,b\rangle^-$ for all $a\in A$ and $b\in B$. In fact, this is how the involution is defined on the dual $\widehat A$. \oldcomment{We have used $\lambda^-$
before. Do we need to mention this again?}
\nl
\bf The $^*$-representation of the Heisenberg algebra  \rm
\nl
Fix a positive right integral $\psi$ on $A$. It is faithful and so $\psi(a^*a)=0$ is only possible when $a=0$.

\notat
Define a scalar product $\langle\,\cdot\, , \,\cdot\, \rangle$ on $A$  by $\langle a',a\rangle=\psi(a^*a')$ for $a,a'\in A$. Denote by $\mathcal H$ the Hilbert space completion of $A$ with respect to this scalar product and let $\Lambda$ be the canonical injection of $A$ in $\mathcal H$. 
\enotat

We have a \emph{problem with the notations}. Here we use $\langle\,\cdot\,,\,\cdot\,\rangle$ for the scalar product while before we have always used this for the pairing between $A$ and $B$. We will keep using the same notation for the two cases because it is so natural. It should be clear from the context which of the two we are dealing with.
\ssnl
Now we consider the GNS-representation of $A$ associated with $\psi$ as in the following proposition. 

\prop\label{prop:repA}
There is a non-degenerate $^*$-representation $\pi$ of $A$ by {\it bounded} operators on $\mathcal H$, given by $\pi(a)\Lambda(x)=\Lambda(ax)$ for $a,x\in A$.
\eprop

\bew
i) Because the right integral is faithful, the canonical map $\Lambda:A\to\mathcal H$ is  injective. Then $\pi(a)$ is well-defined on $\Lambda(A)$. 
\ssnl
ii) We claim that $\pi(a):\Lambda(x)\mapsto \Lambda(ax)$ is a bounded operator. To prove this, 
fix $a \in A$ and consider the map 
\begin{equation}
\Lambda(x)\mapsto \sum_{(a)}\Lambda(a_{(1)})\ot \Lambda(a_{(2)}x)\label{eqn:map}
\end{equation}
from $\Lambda(A)$ to $\mathcal H\ot \mathcal H$. Remark that $\Delta(a)(1\ot x)$ belongs to $A\ot A$ so that this map is well-defined. We have
\begin{align*}\|\sum_{(a)}\Lambda(a_{(1)})\ot \Lambda(a_{(2)}x)\|^2
&=(\psi\ot\psi)((1\ot x^*)\Delta(a)^*\Delta(a)(1\ot x))\\
&=\psi(a^*a)\psi(x^*x)=\psi(a^*a)\|x\|^2.
\end{align*}
We have used that $\psi$ is right invariant. This proves that the map in (\ref{eqn:map}) is continuous. If we take the scalar product in the first factor with $\Lambda(a')$ we will get that also the map 
$\Lambda(x)\mapsto \sum_{(a)}\psi({a'}^*a_{(1)})\Lambda(a_{(2)}x)$
is continuous for all $a,a'\in A$. Because $(A\ot 1)\Delta(A)=A\ot A$ we have that 
any element in $A$ is the linear span of elements of the form
\begin{equation*}
\sum_{(a)}\psi({a'}^*a_{(1)})a_{(2)}.
\end{equation*}
It follows that $\Lambda(x)\mapsto \Lambda(ax)$ is continuous for all $a\in A$. 
\ssnl
iii) We have a unique extension to a bounded operator on all of $\mathcal H$. We also use $\pi(a)$ for this extension on all of $\mathcal H$.  
\ssnl
iv) We clearly get a representation of $A$ by bounded operators. Moreover, by definition we have 
$$\langle\Lambda(ax'),\Lambda(x)\rangle=\psi(x^*ax')=\psi((a^*x)^*x')=
\langle\Lambda(x'),\Lambda(a^*x)\rangle$$
whenever $a,x,x'\in A$. Therefore we have a $^*$-representation. 
\ssnl
v) The representation is non-degenerate because $A^2=A$.
\ebew

The reader should be aware of the fact that the above result is not automatic. One can have faithful positive linear functionals on a $^*$-algebra so that the associated GNS-representation does not provide bounded operators.
\ssnl
In the next proposition, we show that also the dual $^*$-algebra acts on the Hilbert space by means of bounded operators. This result is even less obvious. 

\prop\label{prop:repB}
 There is a non-degenerate $^*$-representation $\gamma$ of the dual algebra $B$ by bounded operators on $\mathcal H$ given by 
\begin{equation*}
\gamma(b)\Lambda(x)=\Lambda(b\tr x)=\sum_{(x)}\langle x_{(2)},b\rangle \Lambda(x_{(1)})
\end{equation*}
where $x\in A$ and $b\in B$.
\eprop

\bew
i) First remark that in the definition of $\gamma(b)$ we use the original pairing of $A$ with its dual $B$. Also the use of the Sweedler notation is justified as the element $x_{(2)}$ is covered by $b$ via the pairing. So $\gamma(b)\Lambda(x)$ is well-defined.
\ssnl
ii) To prove that $\gamma(b)$, defined on $\Lambda(A)$, is a bounded operator we proceed in a similar way as in the proof of  the previous proposition. We take  $a\in A$ and now consider the map
\begin{equation}
\Lambda(x)\mapsto \sum_{(x)}\Lambda(x_{(1)})\ot \Lambda(x_{(2)}a)\label{eqn:map2}
\end{equation}
from $\Lambda(A)$ to $\mathcal H\ot \mathcal H$. Here we use  that $\Delta(x)(1\ot a)$ belongs to $A\ot A$. As before we get 
\begin{equation*}
\|\sum_{(x)}\Lambda(x_{(1)})\ot \Lambda(x_{(2)}a)\|^2
=\psi(a^*a)\psi(x^*x)=\psi(a^*a)\|x\|^2.
\end{equation*}
This shows that also the map in (\ref{eqn:map2}) is continuous. Here we take the scalar product with $\Lambda(a')$ in the second factor. We find that the map
\begin{equation*}
\Lambda(x)\mapsto \sum_{(x)}\psi({a'}^*x_{(2)}a)\Lambda(x_{(1)})
\end{equation*}
is continuous for all $a,a'\in A$ as well. We now claim that linear functionals of the form $\psi({a'}^*\,\cdot\,a)$ span all of the dual $\widehat A$. Indeed, 
\begin{equation*}
\psi({a'}^*\,\cdot\,a)=\psi(\,\cdot\,a\sigma'({a'}^*))
\end{equation*}
where $\sigma'$ is the modular automorphism of $\psi$. Because $A^2=A$, the claim follows. This implies the continuity of the map $\Lambda(x)\mapsto \Lambda(b\tr x)$ for all $b\in B$. We can extend it and we clearly get a bounded representation of $B$ on $\mathcal H$.
\ssnl
iii) We show that it is a $^*$-representation. Therefore take $x,x'\in A$ and $b\in B$. Then
\begin{align*}
\langle \gamma(b)\Lambda(x'),\Lambda(x)\rangle
&=\sum_{(x')}\langle x'_{(2)},b\rangle \langle\Lambda(x'_{(1)}),\Lambda(x)\rangle\\
&=\sum_{(x')}\langle x'_{(2)},b\rangle \psi(x^*x'_{(1)})\\
&=\sum_{(x')}\langle S^{-1}(x^*_{(2)}),b\rangle \psi(x^*_{(1)}x')\\
&=\sum_{(x')}\langle x_{(2)},b^*\rangle^- \psi(x^*_{(1)}x')\\
&=\langle \Lambda(x'),\gamma(b^*)\Lambda(x)\rangle.
\end{align*}
In this series of equalities, we use the brackets for the pairing, as well as for the scalar product. 
We also have used that 
$$(\psi\ot\iota)((a'\ot 1)\Delta(a))=S^{-1}((\psi\ot\iota)(\Delta(a')(a\ot 1)))$$
 for all $a,a'\in A$, see Equation (1.3) in Proposition 1.1 of \cite{VD-part1}. 
\ssnl
iv) Finally, because the action of $B$ on $A$ is unital, so that $A=B\tr A$, we have that $\gamma(B)\mathcal H$ is dense in $\mathcal H$. So the $^*$-representation $\pi$ of $B$ on $\mathcal H$ is non-degenerate.
\ebew

We see from the proof of item iii) that the formula in Equation 1.3 of Proposition 1.1 of \cite{VD-part1} we have used, is essentially the same as saying the $\gamma$ is a $^*$-representation. It justifies the definition of the involution on the dual space.
\ssnl
That these two algebras act on the Hilbert space by means of bounded operators should not come as a surprise. After all, they are the legs of a unitary operator $V$, defined in the  next proposition.

\prop\label{prop:2.5}
The map linear map $V$, defined on $\Lambda(A)\ot\Lambda(A)$ by
\begin{equation*}
V:\Lambda(x)\ot \Lambda(x')\mapsto \sum_{(x)} \Lambda(x_{(1)})\ot \Lambda(x_{(2)}x')
\end{equation*}
is isometric. It extends to a unitary operator on the Hilbert space tensor product $\mathcal H\ot \mathcal H$.
\eprop

\bew
i) As in the proofs of the previous propositions, we have 
\begin{align*}
\| \sum_{(x)} \Lambda(x_{(1)})\ot \Lambda(x_{(2)}x') \|^2
&=(\psi\ot\psi)((1\ot {x'}^*) \Delta(x^*x)(1\ot x'))\\
&=\psi({x'}^*x')\psi(x^*x)
\end{align*}
and we see that $\| V(\Lambda(x)\ot \Lambda(x'))\|=\| \Lambda(x)\ot \Lambda(x')\|$. This is not sufficient to prove that $V$ is isometric, but it is easy to see one can take sums of such elements and the argument will be essentially the same.
\ssnl
ii) Because $\Delta(A)(1\ot A)=A\ot A$, we see that the $V(\Lambda(A)\ot \Lambda(A))$ is dense in $\mathcal H\ot \mathcal H$. Hence $V$ extends to a unitary on this Hilbert space tensor product.
\ebew

The arguments in the proofs of Propositions \ref{prop:repA} and \ref{prop:repB} are implicitly using this.
\ssnl
Again we have a notational ambiguity. 
In \cite{VD-part1}, we have used $V$ for the dualtiy and in particular, in Proposition 3.17 of \cite{VD-part1}, 
we have shown however that it acts on $A\ot A$ as the canonical map $a\ot a'\mapsto \Delta(a)(1\ot a')$. In the previous proposition, we consider this map on the Hilbert space level. This justifies the use of the same symbol. This remark is related to the next one.

\notat
i) Both the representation $\pi$ of $A$ and $\gamma$ of $B$ on the Hilbert space $\mathcal H$ are faithful. Therefore, it makes sense to omit the use of  $\pi$ for $A$ and $\gamma$ for $B$ and let these algebras act directly on $\mathcal H$. So we write $a\xi$ and $b\xi$ for $\pi(a)\xi$ and $\gamma(b)\xi$ when $a\in A$ and $b\in B$ and where $\xi$ is any vector in $\mathcal H$. Of course in the case where $\xi=\Lambda(x)$, with $x\in A$ we have
\begin{equation*}
a\Lambda(x)=\Lambda(ax)
\tussenen
b\Lambda(x)=\sum_{(x)}\langle x_{(2)},b\rangle \Lambda(x_{(1)})
\end{equation*}
for $a\in A$ and $b\in B$.
\ssnl
ii) With this convention, the operator $V$ on $\Lambda(A)\ot\Lambda(A)$ is just the same as in the canonical map $T$, as in Proposition 3.17 in \cite{VD-part1}.
\enotat
 
The two $^*$-representations $\pi$ of $A$ and $\gamma$ of $B$ combine to a non-degenerate $^*$-representation of the Heisenberg algebra $AB$ on $\mathcal H$ with bounded operators. 
\ssnl
In Proposition 3.16 of \cite{VD-part1}, we have seen that the Heisenberg algebra acts on $A$ by means of operators of the form $x\mapsto \langle x,b\rangle a$ 
 where $a\in A$ and $b\in B$. Because any element in $B$ is of the form $\psi({a'}^*\,\cdot\,)$, in its representation on the Hilbert space, we get operators of the form $\xi\mapsto \langle \xi,\Lambda(a')\rangle\Lambda(a)$ with $a,a'\in A$. The operators are bounded rank one operators. Observe that not all rank one operators are of this type because not all vectors in $\mathcal H$ belong to the range of $\Lambda$. However, because this range is dense, we get the following result.

\prop
 The norm closure of the Heisenberg algebra in the representation is the C$^*$-algebra of all compact operators and the weak closure is the von Neumann algebra of all bounded operators on the Hilbert space $\mathcal H$.
 \eprop
\keepcomment
{Moet hier nog iets meer bij? Referenties? Intern/ extern? Nog iets opmerken over regularity of the multiplicative unitary?}
\nl
\bf The coproduct on the operator algebra completions of $A$ \rm
\nl
In this item, we use the unitary operator $V$ to get extensions of the coproduct to the following operator algebras.

\notat\label{notat:2.8}
We denote by $P$ the {\it norm closure} of $A$ as acting on the Hilbert space $\mathcal H$ and we use $M$ for the {\it weak closure} of $A$. In other words, $P$ is the C$^*$-algebra and $M$ the von Neumann algebra generated by $A$ on $\mathcal H$. 
\enotat

In Section \ref{sect:dual} we will consider the operator algebras associated to the dual, generated by $B$ as acting on $\mathcal H$, see Notation \ref{notat:4.5}.
\ssnl
The coproduct $\Delta$ on $A$ has an obvious extension to the C$^*$-algebra $P$ and the von Neumann algebra $M$. Before we formulate this result, we recall the two notions of a coproduct in the operator algebraic framework.

\defin\label{defin:2.9}
i) Let $P$ be any C$^*$-algebra. A coproduct on $P$ is a non-degenerate $^*$-homomorphism $\Delta$ from $P$ to the multiplier algebra 
$M(P\ot P)$ of the minimal C$^*$-tensor product $P\ot P$. It is assumed to be coassociative. The sets $\Delta(P)(1\ot P)$ and $(P\ot 1)\Delta(P)$ are supposed to be contained in $P\ot P$.
\ssnl
ii) If $M$ is any von Neumann algebra, a coproduct on $M$ is a unital and normal $^*$-homomorphism $\Delta:M\to M\ot M$ where now the von Neumann algebra tensor product is considered. Again it is assumed to be coassociative.
\edefin

These concepts are found in various places in the literature on locally compact quantum groups. See e.g.\ \cite{K-V1}, \cite{K-V2} and \cite{K-V3}. In the more recent work (see \cite{VD-sigma} and \cite{VD-notes}), there is included some discussion about the two notions.
\oldcomment{
\rood\ssnl
Deze opmerkingen moeten we misschien wel eens verifiëren.
\ssnl}{}
\ssnl
With these definitions, we can now formulate the following result.  Here again $P$ is the C$^*$-algebra generated by $A$  and $M$ the von Neumann algebra generated by $A$, both as acting on the Hilbert space. 

\prop\label{prop:2.10a}
Define $\Delta_1$ on $M$ by 
$$\Delta_1(x)=V(x\ot 1)V^*.$$ 
Then $\Delta_1$ is a coproduct on the von Neumann algebra $M$ (as in ii) of Definition \ref{defin:2.9}).
The restriction of $\Delta_1$ to  $P$ is a coproduct on the C$^*$-algebra (as in i) of Definition \ref{defin:2.9}). It extends the original coproduct on $A$.
\eprop

\bew
i) If $a$ is an element in $A$, we know that $V(a \ot 1)V^*=\Delta(a)$ in the multiplier algebra of $AB\ot A$. The right hand side belongs to the algebraic multiplier algebra $M(A\ot A)$. We clearly have that $V(a\ot 1)V^*$ is well-defined and bounded on all of $\mathcal H\ot\mathcal H$. On the other hand, $\Delta(a)$ is well-defined on vectors of the form $\xi\ot \Lambda(x)$ because $\Delta(a)(1\ot x)\in A\ot A$ for all $a,x\in A$. It is again a bounded operator and extends to all of $\mathcal H\ot \mathcal H$.
\ssnl
 \oldcomment{Argue that $\Delta(a)$ acts as a bounded operator. Refer to here in Section 4. See a remark after Proposition \ref{prop:2.11}. } 

We now view this equation as bounded operators on $\mathcal H\ot \mathcal H$. Then $\Delta(a)$ belongs to the multiplier algebra $M(P\ot P)$ of the C$^*$-algebra tensor product for all $a\in A$. 
\ssnl
 Furthermore, since the C$^*$-algebra $P$ acts in a non-degenerate way, this multiplier algebra is a subalgebra of the von Neumann tensor product $M\ot M$. By continuity, also $V(x\ot 1)V^*$ will belong to $M\ot M$ for all $x\in M$. We clearly have that $\Delta_1$ is a normal, unital $^*$-homomorphism from $M$ to $M\ot M$. Coassociativity on $M$  follows by continuity from the coassociativity of $\Delta$ on $A$. It is also a consequence of the fact that $V$ satisfies the Pentagon equation. This proves the first statement of the proposition.
\vskip 3pt
ii) We have just argued that $\Delta_1(x)\in M(P\ot P)$ for $x\in A$ and again by continuity, we have that this is still true for $x\in P$. Because $\Delta(A)(1\ot A)$ and $(A\ot 1)\Delta(A)$ are equal to the algebraic tensor product $A\ot A$, we will have that  $\Delta_1(P)(1\ot P)$ and $(P\ot 1)\Delta_1(P)$ will be subsets of the C$^*$-tensor product $P\ot P$. This proves that we have a genuine coproduct on the C$^*$-algebra $P$.
\ebew

\keepcomment{We have to be more careful here. Check all statements and write them in detail for ourselves. \rood To Do!
\ssnl}
In what follows, we will also use $\Delta$ to denote these extensions of the original coproduct on $A$ to $M$, and so also to $P$.
\ssnl
In Section \ref{sect:dual} we will do the same for the dual $B$. 
\ssnl
First we extend the integrals to Haar weights on the von Neumann algebra $M$. This is done in the next section.

\section{\hspace{-17pt}. The Haar weights on the pair $(M,\Delta)$}\label{sect:hwt} % \input artikel2b.tex %\newpage

In the previous section we have extended the coproduct from $A$ to the von Neumann algebra $M$, generated by $A$ in the GNS-representation of the right integral $\psi$ on $A$. In this section we extend $\psi$ to a faithful normal semi-finite weight on $M$ and show that it is a right Haar weight. We use more properties of the original algebraic quantum group to obtain the left Haar weight. The result is a locally compact quantum group.
\ssnl 
The extension of the integrals to Haar weights on the C$^*$-algebra and the von Neumann algebra completions is obtained by using standard techniques from the theory of left Hilbert algebras. This is not so simple but unfortunately, it seems not possible to avoid this theory. 
\ssnl
The material is found in Paragraph 1 in Chapter VI  on left and right Hilbert algebras in \cite{T3}. As we do not need the full strength of this theory, and for the convenience of the reader who is not familiar with it, we have included two appendices where we discuss that part that is needed to get the further results.  We will mostly refer to these appendices instead of to \cite{T3}. References to \cite{T3} are found in the appendices.
\nl
It is easy to prove the following result. For the definition of a left Hilbert algebra we refer to Definition \ref{defin:lha} in  Appendix \ref{sect:appA}.

\prop\label{prop:lha}
The space $\Lambda(A)$ is a left Hilbert algebra $\mathfrak A$ with the $^*$-algebra structure inherited from $A$ and with the given scalar product. 
\eprop

\bew
Most of the axioms in Definition \ref{defin:lha} are trivial reformulations of earlier results on the GNS-representation of $A$ on the Hilbert space $\mathcal H$ defined as the completion of $\Lambda(A)$. The density of $\mathfrak A^2$ in $\mathcal H$ is a consequence of the fact $A^2=A$. 
\ssnl
Essentially we only have to prove that the $^\sharp$-operation is preclosed. For this we take elements $a,a'\in A$ and we write
\begin{align}
\langle \Lambda(a^*),\Lambda(a')\rangle
&=\psi({a'}^*a^*)=\psi(a^*\sigma'({a'}^*))\\
&=\langle\Lambda(\sigma'({a'}^*)),\Lambda(a)\rangle.
\end{align}
We have used the modular automorphism $\sigma'$ of $\psi$. This proves that $\Lambda(a')$ is in the domain $\mathcal D^\flat$ of the adjoint $^\flat$ of the $^\sharp$-map and that 
\begin{equation*}
 (\Lambda(a'))^\flat=\Lambda(\sigma'({a'}^*)).
\end{equation*}
This completes the proof.
\ebew

Recall that we are working here with conjugate linear operators. For the adjoint $^\flat$ of $^\sharp$ we have the basic equatity
\begin{equation*}
\langle \xi_1^\sharp, \xi_2\rangle=\langle\xi_1,\xi^\flat_2\rangle^-=\langle\xi_2^\flat, \xi_1\rangle
\end{equation*}
when $\xi_1\in \mathcal D^\sharp$ and $\xi_2\in \mathcal D^\flat$.
\ssnl
We have used the same symbol for the Hilbert space completion of the left Hilbert algebra because this clearly is the same as the Hilbert space completion of $\Lambda(A)$ as used before.
\ssnl
Remark that it is expected, as we see in the proof above, that the existence of the modular automorphism on the algebra $A$ for the positive linear functional $\psi$ is sufficient to get that the $^\sharp$-operator has a densely defined adjoint. 
\nl
For the construction of the Haar weights we use the procedure as described in the appendix. 
\ssnl
First we  need right bounded vectors. Recall the following definition in the framework of left Hilbert algebras (cf. Definition \ref{defin:rbe} in Appendix \ref{sect:appA}).

\defin
A vector $\eta\in \mathcal H$ is called {\it right bounded} if there is a bounded operator on $\mathcal H$, necessarily unique and denoted by $\pi_r(\eta)$, satisfying $\pi_r(\eta)\xi=\pi_\ell(\xi)\eta$ for all $\xi\in \mathfrak A$.
\edefin

For the left Hilbert algebra we consider here, a right bounded vector is an element $\eta$ of $\mathcal H$ such that the map $\Lambda(a)\mapsto a\eta$ is bounded. The unique bounded extension is the operator $\pi_r(\eta)$ and so we have $\pi_r(\eta)\Lambda(a)=a\eta$ for all $a\in A$.
\ssnl
The integral $\psi$ can now be extended to a normal semi-finite faithful weight on the von Neumann algebra completion $M$ of $A$. This is an application of the main result, Theorem \ref{stel:main} in  Appendix \ref{sect:appA}.

\stel\label{stel:extw}
Define $\psi_1(x)$ for $x\in M$ and $x\geq 0$ by
\begin{equation}
\psi_1(x)=\sup\{\langle x\eta,\eta\rangle \mid \eta \text{ is right bounded and }\|\pi_r(\eta)\|\leq 1\}.\label{eqn:fnsw}
\end{equation}
Then $\psi_1$ is a faithful normal semi-finite weight on the von Neumann algebra $M$.  If $a\in A$, then $\psi_1(a^*a)=\psi(a^*a)$.
\estel

\bew
The result is a straightforward consequence of the general theorem in the appendix. The last statement follows because we have $\psi_1(x^*x)=\langle\xi,\xi\rangle$ when $\xi\in\mathfrak A$ and $x=\pi_\ell(\xi)$. Indeed, we have $a=\pi_\ell(\xi)$ when $\xi=\Lambda(a)$ where $a\in A$.
\ebew

\opm
One of the problems to show that Equation {\ref{eqn:fnsw}}) actually defines a weight, is to prove  that it is additive. This is based on Lemma \ref{lem:A.21} of Appendix \ref{sect:appB} where this result is proven. One may wonder if it would be easier to obtain the additivity  from the knowledge that it is additive on the dense subset of positive elements in $A$. But apparently, this does not seem to be the case.
\oldcomment{The additivity can probably also be proven by approximation. We have made a remark about this in the appendix. This would make the argument simpler and avoid a difficult part in the proof there. \rood To investigate! }{}
\eopm

In what follows, \emph{we will use $\psi$ also for the weight $\psi_1$ on the von Neumann algebra}.
\ssnl
For the following result, we need left bounded vectors (see Definition \ref{defin:A.16a} in Appendix \ref{sect:appA}).

\defin
An element $\xi$ in $\mathcal H$ is called \emph{left bounded} if there is a bounded operator, necessarily unique and denoted by $\pi_\ell(\xi)$ satisfying $\pi_\ell(\xi)\eta=\pi_r(\eta)\xi$ for all right bounded elements $\eta$.
\edefin

Any element  $\xi\in\mathfrak A$ is left bounded because we have $\pi_\ell(\xi)\eta=\pi_r(\eta)\xi$ by the definition of right bounded elements. We also see that there is no conflict with the notations. 

\prop\label{prop:2.16}
Denote by $\mathfrak N_\psi$ the set of elements $x\in M$ satisfying $\psi(x^*x)\leq \infty$. They are precisely the operators $\pi_\ell(\xi)$ where $\xi$ is a \emph{left bounded vector}. 
\eprop

This is Proposition \ref{prop:A.23a} in  Appendix \ref{sect:appA}. 
\ssnl
It is a also consequence of this result that the canonical Hilbert space $\mathcal H_\psi$ is the same as the original Hilbert space $\mathcal H$ and that the canonical injection $\Lambda_\psi$ from the left ideal $\mathfrak N_\psi$ to $\mathcal H_\psi$ is given by $\pi_\ell(\xi)\mapsto \xi$. See also Proposition \ref{prop:A.21}.
\ssnl
 We will need this to show that the weight $\psi$ on the von Neumann algebra is still right invariant in the following proposition. \oldcomment{Formulate as a result.}

\oldcomment{\rood Formulate/ recall results from the appendix, needed to prove the following.}{}

\prop\label{prop:2.18}
The weight $\psi$ as defined in Theorem \ref{stel:extw} is right invariant on $(M,\Delta)$. 
\eprop

\bew
i) Let $a,a'$ be elements in $A$.  Take a right bounded vector $\eta$ in $\mathcal H$. Then we have
\begin{align*}
\Delta(a)(\eta\ot\Lambda(a'))
&=\sum_{(a)}a_{(1)}\eta\ot \Lambda(a_{(2)}a')\\
&=\sum_{(a)}\pi_r(\eta)\Lambda(a_{(1)})\ot \Lambda(a_{(2)}a')\\
&=(\pi_r(\eta)\ot 1)V(\Lambda(a)\ot\Lambda(a')).
\end{align*}
By continuity and the density of $\Lambda(A)$ in $\mathcal H$ we find
\begin{equation}
\Delta(a)(\eta\ot\xi)=(\pi_r(\eta)\ot 1)V(\Lambda(a)\ot\xi)\label{eqn:2.6}
\end{equation}
for all $\xi\in\mathcal H$.
\ssnl
ii) Now we use a property of the theory of left Hilbert algebras. Consider $x\in\mathfrak N_\psi$. As we have seen in Proposition \ref{prop:2.16}, there is a left bounded vector $\xi$ such that $\pi_\ell(\xi)=x$. By Proposition \ref{prop:A.23} in Appendix \ref{sect:appA}, there is a bounded sequence $(a_n)$ in $A$ such that $\Lambda(a_n)\to\Lambda(x)$.
\oldcomment{\rood Formulate the result in the appendix like this.} For all right bounded vectors $\eta$ we have
\begin{equation*}
 a_n\eta=\pi_r(\eta)\Lambda(a_n)\to\pi_r(\eta)\Lambda(x)=x\eta. 
\end{equation*}
Therefore $(a_n\ot 1)\xi_1\to (x\ot 1)\xi_1$ for a dense set of vectors $\xi_1$ in $\mathcal H\ot \mathcal H$. Because the sequence $(a_n\ot 1)$ is bounded, we have this for all vectors $\xi_1$. Therefore
\begin{equation*}
(a_n\ot 1)V^*(\eta\ot \xi)\to (x\ot 1)V^*(\eta\ot \xi)
\end{equation*}
for all vectors $\eta$ and $\xi$.
If we apply $V$ we see that $\Delta(a_n)(\eta\ot\xi)\to \Delta(x)(\eta\ot \xi)$. We now use Equation (\ref{eqn:2.6}) for $a_n$ and we take the limit. Because also $\Lambda(a_n)\to\Lambda(x)$ It follows that still
\begin{equation*}
\Delta(x)(\eta\ot\xi)=(\pi_r(\eta)\ot 1)V(\Lambda(x)\ot\xi)
\end{equation*}
for all right bounded vectors $\eta$, all $\xi\in\mathcal H$ and all $x\in  \mathfrak N_\psi$.
\ssnl
iii) Then we obtain, with $\omega=\langle\,\cdot\,\xi,\xi\rangle$, and the supremum taken over right bounded vectors $\eta$ satisfying $\|\pi_r(\eta)\|\leq 1$, 
\begin{align*}
\psi((\iota\ot\omega)\Delta(x^*x))
&=\sup \langle \Delta(x^*x)(\eta\ot\xi),\eta\ot \xi\rangle\\
&=\sup \langle (\pi_r(\eta)\ot 1)V(\Lambda(x)\ot\xi),(\pi_r(\eta)\ot 1)V(\Lambda(x)\ot\xi)\rangle\\
&= \langle V(\Lambda(x)\ot\xi),V(\Lambda(x)\ot\xi)\rangle\\
&=\langle \Lambda(x)\ot\xi,\Lambda(x)\ot\xi\rangle\\
&=\psi(x^*x)\omega(1)
\end{align*}
This is sufficient to have that $\psi$ is invariant. 
\keepcomment{\rood Nog wat uitleg is nodig.}{}
\ebew

This gives the right Haar weight on the pair $(M,\Delta)$. We now look for the left Haar weight.
\nl
\bf The left Haar weight on $(M,\Delta)$ \rm
\nl
We now need to extend the left integral $\varphi$, defined on $A$, to the von Neumann algebra. There are several ways to do this. 
\ssnl 
One approach would be the following. We could apply the previous arguments to the algebraic quantum group obtained from $A$ by flipping the coproduct. The original left integral becomes the right integral. However, it is not clear that the resulting von Neumann algebra completion is the same as the one obtained originally from the right integral. It requires some work to show this. 
\keepcomment{We refer to our lecture notes for a situation where this method is used .  See Proposition 2.24 op page 52 in the version of May 2016.
\ssnl
Aan te vullen, eventueel meer uit te leggen. Waarom dit nuttig is, ook al lijkt er op voorhand wat overbodig werk bij te komen?}
\ssnl
A second possibility is to use methods from the general theory of locally compact quantum groups. Doing so, we need the operator algebraic version of the antipode and related objects. This method has the advantage that we obtain other important properties of the locally compact quantum group $(M,\Delta)$.  \oldcomment{\rood We will say more about this later in this section.}{}
\ssnl
However, completely in line of the scope of this paper, we will make optimal use of  the properties of the original algebraic quantum group $(A,\Delta)$ as we have reviewed in Section \ref{sect:prel}. In particular, the square root $\delta^\frac12$ of the modular element $\delta$ obtained in Proposition \ref{prop:1.11a}, is what we use for this method. 
\ssnl
Again we use the theory of left Hilbert algebras.

\prop 
Define a map $\Lambda_\varphi$ from $A$ to $\mathcal H$ by $\Lambda_\varphi(a)=\Lambda(a\delta^\frac12)$. Then $\Lambda_\varphi(A)$ is a left Hilbert algebra  for the $^*$- algebra structure inherited from $A$ and with the given scalar product form $\mathcal H$. The left von Neumann algebra of this left Hilbert algebra is $M$.
\eprop

\bew
The proof is mostly straightforward, just as it was for the case earlier in Proposition \ref{prop:lha}.
\ssnl
First remark that $\delta^\frac12$ is a well-defined multiplier of $A$ as was shown in Proposition \ref{prop:1.11a}. If $\xi,\eta$ are elements in this left Hilbert algebra, given by $\xi=\Lambda_\varphi(a)$ and $\eta=\Lambda_\varphi(c)$, then 
\begin{equation*}
\xi\eta=\Lambda_\varphi(ac)=\Lambda(ac\delta^\frac12)=a\Lambda(c\delta^\frac12)=a\Lambda_\varphi(c)=a\eta
\end{equation*}
and we see that $\pi_\ell(\xi)=a$ as acting on the Hilbert space.
\ssnl
For the adjoint of the $^\sharp$ map we calculate, for all $a,c\in A$, 
\begin{align*}
\langle\Lambda_\varphi(a^*),\Lambda(c)\rangle
&=\langle\Lambda(a^*\delta^\frac12),\Lambda(c)\rangle
=\psi (c^*a^*\delta^\frac12)\\
&=\psi(a^*\delta^\frac12\sigma'(c^*))
=\psi (\delta^\frac12 a^*\delta^\frac12\sigma'(c^*)\delta^{-\frac12})\\
&=\langle\Lambda(\delta^\frac12\sigma'(c^*)\delta^{-\frac12}) ,\Lambda(a\delta^\frac12)\rangle\\
&=\langle\Lambda(\delta^\frac12\sigma'(c^*)\delta^{-\frac12}) ,\Lambda_\varphi(a)\rangle.
\end{align*}
We have used that $\delta^\frac12$ is self-adjoint and that $\sigma'(\delta^\frac12)=\delta^\frac12$.
\ssnl
This proves that $\Lambda(A)$ belongs to the domain of the adjoint of the $^\sharp$-map.
\ebew

Next consider the Haar weight $\varphi_1$ on $M$ obtained for this left Hilbert algebra.

\prop
For all $a\in A$ we have $\varphi_1(a^*a)=\varphi(a^*a)$. It is a left invariant weight on $(M,\Delta)$.
\eprop

We clearly have, for all $a\in A$. 
\begin{align*}
\varphi_1(a^*a)
&=\langle\Lambda_\varphi(a),\Lambda_\varphi(a)\rangle\\
&=\langle\Lambda(a\delta^\frac12),\Lambda(a\delta^\frac12)\rangle\\
&=\psi(\delta^\frac12 a^*a \delta^\frac12)\\
&=\psi(a^*a\delta)=\varphi(a^*a).
\end{align*}

The left invariance of $\varphi_1$ on $M$ is obtained from the left invariance of $\varphi$ on $A$, just as in the proof of Proposition \ref{prop:2.18}. \keepcomment{We can refer to the temporary file tmp10 for details for ourselves.}{}
\ssnl
In what follows we use also $\varphi$ for the faithful normal semi-finite weight $\varphi_1$ on $M$, just as we did for the right Haar weight.
Now we are ready for the main result in this section.
\nl
\bf The main result \rm
\nl
We arrive at the following.
\stel
The pair $(M,\Delta)$ is a locally compact quantum group with right Haar weight $\psi$ and left Haar weight $\varphi$.
%=\psi\circ R$.
\estel

There is now the general procedure to pass from the von Neumann algebraic quantum group $(M,\Delta)$ to a C$^*$-algebraic one. In this case we get what is expected.

\stel
The C$^*$-algebraic quantum group obtained from $(M,\Delta)$ using the general theory is the pair $(P,\Delta)$ where $P$ is the C$^*$-algebra completion of $A$ as considered before and where $\Delta$ is the restriction of the coproduct on $M$ to $P$.
\estel

\bew
From the general theory, we know that the associated C$^*$-algebra is obtained as the closure of the space of elements $(\omega\ot\iota)V$ where $\omega$ runs through all normal linear functionals on $\mathcal B(\mathcal H)$. We can restrict to functionals of the type $\omega=\langle \,\cdot \,\Lambda(a),\Lambda(b)\rangle$ where $a,b\in A$. 
Because $V(\Lambda(a)\ot \xi)=\sum_{(a)}\Lambda(a_{(1)}\ot a_{(2)}\xi$ then we get
\begin{equation*}
(\omega\ot\iota)V=\sum_{(a)}\langle \Lambda(a_{(1)}),\Lambda(b)\rangle a_{(2)} =\sum_{(a)}\psi(b^*a_{(1)})a_{(2)}.
\end{equation*}
Such elements span all of $A$ and we see that the associated $C^*$-algebra indeed coincides with the closure $P$ of $A$.
\ebew

\oldcomment{Is this ok?}{\ssnl}

It is interesting to investigate the further results, obtained from the general theory of locally compact quantum groups, and how they are formulated in terms of the original algebraic quantum group. One of the issues to consider is the analytic structure of the algebraic quantum group, as obtained in Section \ref{sect:prel} in view of the properties of the locally compact quantum group. These results are reviewed in the third part of this series \cite{VD-part3}.
 \ssnl
Here we now will focus on the construction of the dual locally compact quantum group. This is the scope of the next section.

\section{\hspace{-17pt}. The dual and the Haar weights on the dual}\label{sect:dual} %\input artikel2c.tex %\newpage

To study the dual  of the locally compact quantum group $(M,\Delta)$ that we have obtained in the previous two sections, we can choose between \emph{two methods}. We can either first consider the dual of the algebraic quantum group and then apply the completion procedure as in the two preceding sections to obtain the locally compact quantum group completion. Or we can simply take the dual of the completed locally compact quantum group $(M,\Delta)$ as in the general theory. 
\ssnl
What we essentially do in this section is to explain that these two methods provide the same result. We will follow the first method and explain at the various levels that we get the same result when we apply the second method.
\ssnl
Before we start this,  we have to mention some of the basic differences between the treatment of the dual in  the general theory of locally compact quantum groups and the approach here.

\opm
i) In the general theory (as e.g.\ explained in \cite{VD-sigma}), the left Haar weight is used to begin with. Here we started with  the right integral and thus obtained the right Haar weight on the von Neumann algebra.
\ssnl
ii) For the coproduct on the dual, in the theory of locally compact quantum groups, one usually takes the opposite coproduct of the one taken in the algebraic theories.
\ssnl
iii) The left Haar weight on the dual is constructed first.
\eopm

We will stick to the algebraic convention because that is common in the theory of algebraic quantum groups. This implies however that we have to modify the formulas when we refer to the duality as obtained in the theory of locally compact quantum groups. But the changes are quite obvious and this should not present any real difficulty.
\nl
As before, we let $(B,\Delta)$ be the dual of $(A,\Delta)$ in the sense of algebraic quantum groups. The Fourier transform from $A$ to $B$, defined by $a\mapsto \psi(S(\,\cdot\,)a)$ where $\psi$ is the right integral, plays an important role. For the left integral $\widehat\varphi$ on $B$ we have the Plancerel formula saying that $\widehat\varphi(b^*b)=\psi(a^*a)$ when $b=\psi(S(\,\cdot\,)a)$ for $a\in A$. See item ii) in Proposition \ref{prop:C.10} in Section \ref{sect:prel}.
\ssnl
The first important step in this procedure is motivated by an observation in the theory of locally compact quantum groups (see e.g.\ the remark before Definition 4.6 in Section 4 of \cite{VD-sigma}). It is characteristic for the theory that the spaces $L^2(G)$ and $L^2(\widehat G)$ of square integrable functions on the locally compact abelian group $G$ and its dual $\widehat G$ are identified via the Fourier transform. In this case, it means that we combine the map $a\mapsto \Lambda(a)$ from $A$ to $\mathcal H$ with the Fourier transform to arrive at a map from $B$ to $\mathcal H$. 
\ssnl
This leads us to the following definition and results. 

\defin\label{defin:4.2a}
We define a map $\widehat\Lambda:B\to\mathcal H$ by $\widehat\Lambda(b)=\Lambda(a)$ if $b=\psi(S(\,\cdot\,)a)$ where $a\in A$.
\edefin

For this map we have the following property. We will use it later when we compare the two methods for constructing the left Haar weight on the dual.

\prop
For all $b\in B$ and $x\in A$ we have
\begin{equation*}
\langle\widehat\Lambda(b),\Lambda(x)\rangle=\langle S(x)^*,b\rangle.
\end{equation*}
We use the Hilbert scalar product on the left and the pairing on the right.
\eprop

\bew
Take $a\in A$ and let $b=\psi(S(\,\cdot\,)a)$. Then for all $x\in A$ we have
\begin{align*}
\langle\widehat\Lambda(b),\Lambda(x)\rangle=\langle\Lambda(a),\Lambda(x)\rangle=\psi(x^*a)
\end{align*}
while also
\begin{equation*}
\langle S(x)^*,b\rangle=\psi(S(S(x)^*)a)=\psi(x^*a).
\end{equation*}
\ebew

We have that $\widehat\varphi(b^*b)=\psi(a^*a)$ because of item ii) in Proposition \ref{prop:C.10}. Therefore 
 $\widehat \varphi(b^*b)=\langle\widehat\Lambda(b),\widehat \Lambda(b)\rangle$ for all $b\in B$.
\ssnl
We expect that $\widehat\Lambda$ is the GNS-map associated with the dual left integral $\widehat \varphi$ on the $^*$-algebra $B$. That this is indeed the case and it follows from the following result.

\prop\label{prop:2.29}
For all $b,b'\in B$ we have $\widehat\Lambda(b'b)=b'\widehat\Lambda(b)$.
\eprop

\bew
When $b=\psi(S(\,\cdot\,)a)$ for $a\in A$, we have for all $x\in A$,
\begin{align*}
\langle x,b'b\rangle
&=\sum_{(x)} \langle x_{(1)},b'\rangle \psi(S(x_{(2)})a) \\
&=\sum_{(S(x))} \langle S\inv ((S(x))_{(2)}),b'\rangle \psi((S(x))_{(1)})a) \\
&=\sum_{(a)} \langle  a_{(2)},b'\rangle \psi(S(x)a_{(1)}).
\end{align*}
For the last equality we have used the Equation (1.3) in Proposition 1.1 of \cite{VD-part1}. 
Then 
\begin{align*}
\widehat\Lambda(b'b)
&=\sum_{(a)} \langle  a_{(2)},b'\rangle \Lambda(a_{(1)})\\
&=b'\Lambda(a)=b'\widehat\Lambda(b).
\end{align*}
\vskip-10pt
\ebew

In Notation \ref{notat:2.8} we have introduced the von Neumann algebra $M$ and the C$^*$-subalgebra $P$, respectively the weak and norm closures of $A$ as acting on the Hilbert space $\mathcal H$. We now do the same for $B$.

\notat\label{notat:4.5}
We use $Q$ and $N$ respectively for the C$^*$-algebra and the von Neumann algebra generated by $B$ as acting on $\mathcal H$.
\enotat

In the theory of locally compact quantum groups, the dual von Neumann algebra $\widehat M$ is defined here as the left leg of the unitary operator $V$. This coincides with the von Neumann algebra $N$ as we see from the next result. Similarly for the associated $C^*$-algebras. It is a natural consequence of the fact that $V$ is essentially the duality.

\prop
Let $b=\psi(S(\,\cdot\,)a)$ for $a$ in $A$ and let $\omega$ be the associated linear functional on $M$. Then
$(\iota\ot\omega)V=b$, seen as operators on $\mathcal H$. 
\eprop
\bew
First remark that elements of $B$ yield bounded linear functionals on $M$. Indeed, given $a,a'\in A$ we have 
\begin{align*}
\langle x\Lambda(a'),\Lambda(a) \rangle
&=\psi(a^*xa)=\psi(S(a^*xa)\delta\inv) \\
&=\psi(S(a)S(x)S(a^*)\delta\inv)\\
&=\psi(S(x)S(a^*)\delta\inv\sigma'(S(a)))
\end{align*}
for all $x\in A$.  We have used a formula from Proposition 1.6 in \cite{VD-part1} and that the scaling constant is $1$. Such functionals are continuous and extend to normal functionals on $M$. Because $A^2=A$ in the end we get that all functionals of the form $x\mapsto \psi(S(\,\cdot\,)a)$ are continuous and extend to normal linear functionals on $M$. 
 \keepcomment{Did we use this earlier? I don't think so. Wait for deleting this remark, and check later again!}{}
%\ssnl
For all $x\in A$ and $\xi\in\mathcal H$ we have $V(\Lambda(x)\ot\xi)=\sum_{(x)}\Lambda(x_{(1)})\ot x_{(2)}\xi$ and so
\begin{equation*}
((\iota\ot\omega)V)\Lambda(x)=\sum_{(x)}\omega(x_{(2)})\Lambda(x_{(1)})=\sum_{(x)}\langle x_{(2)},b\rangle \Lambda(x_{(1)})=\gamma(b)\Lambda(x).
\end{equation*}
\vskip-15pt\ebew

Next we show that the coproduct on $B$ is precisely the one we have as in the duality framework of locally compact quantum groups.

\prop\label{prop:4.7a}
For all $b\in B$ and $\xi\in\mathcal H$ of the form $\Lambda(a)$ for $a\in A$ we have 
\begin{equation*}
V^*(\xi\ot\widehat\Lambda(b))=\sum_{(b)}b_{(1)}\xi\ot \widehat\Lambda(b_{(2)}).
\end{equation*}
Consequently we get $V^*(1\ot b)V=\Delta(b)$ as operators on $\mathcal H\ot\mathcal H$.
\eprop

\bew
Let $\xi=\Lambda(a)$ and $b=\psi(S(\,\cdot\,)c)$ with $a,c\in A$. Then
\begin{align*}
V^*(\xi\ot\widehat\Lambda(b))
&=\sum_{(a)}\Lambda(a_{(1)})\ot S(a_{(2)})\Lambda(c) \\
&=\sum_{(a)}\Lambda(a_{(1)})\ot \Lambda(S(a_{(2)})c).
\end{align*}
On the other hand we have 
\begin{equation}
\sum_{(b)}b_{(1)}\Lambda(a)\ot\widehat\Lambda(b_{(2)})
=\sum_{(a)(b)}\Lambda(a_{1)})\ot \langle a_{(2)},b_{(1)}\rangle\widehat\Lambda(b_{(2)}). \label{eqn:2.16}
\end{equation}
Now we have, again for all $a,x$ in $A$ that
\begin{equation*}
\sum_{(b)}\langle a,b_{(1)}\rangle\langle x,b_{(2)}\rangle
=\langle ax,b\rangle=\psi(S(ax)c)=\psi(S(x)S(a)c).
\end{equation*}
Then we find 
\begin{equation*}
\sum_{(b)}\langle a,b_{(1)}\rangle b_{(2)}
=\psi(S(\,\cdot\, )S(a)c)
\end{equation*}
and
\begin{equation*}
\sum_{(b)}\langle a,b_{(1)}\rangle \widehat\Lambda(b_{(2)})
=\Lambda(S(a)c).
\end{equation*}
If we apply this with $a$ replaced by $a_{(2)}$ in Equation (\ref{eqn:2.16}), we find 
\begin{equation*}
V^*(\xi\ot \widehat\Lambda(b))=\sum_{(b)}b_{(1)}\xi\ot\widehat\Lambda(b_{(2)}).
\end{equation*}
\ebew

The use of the Sweedler notation in the above proof requires some extra attention. In particular, in Equation (\ref{eqn:2.16} on has to remark that elements of the form
\begin{equation*}
\sum_{(a)(b)}\langle a_{(2)},b_{(1)}\rangle a_{1)} \ot b_{(2)})
\end{equation*}
are well-defined in $A\ot B$. See e.g.\ the item on the Heisenberg algebra in Section 3 of \cite{VD-part1}.
\ssnl
We also remark that $\Delta(b)$ acts as a bounded operator on $\mathcal H\ot\mathcal H$. This follows from the fact that $B$ acts by bounded operators and because $V^*(1\ot b)V=\Delta(b)$ as acting on $\mathcal H\ot\mathcal H$. See Proposition \ref{prop:2.10a} in Section \ref{sect:lcp} where we have a similar argument for the equation $V(a\ot 1)V^*=\Delta(a)$.
\oldcomment{Verify! Explain We have to make a remark about this. How $\Delta(B)$ acts on $\mathcal H\ot\mathcal H$ by bounded operators. We have to do this also for $\Delta(a)$ before!}{}
\ssnl
Just as in the case of $A$, see Proposition \ref{prop:2.10a}, the coproduct on $B$ can be extended to the C$^*$-algbra $Q$ and the von Neumann algebra $N$ generated by $B$ as represented on the Hilbert space. 

\prop\label{prop:2.11}
Define $\Delta_1$ on $N$ by 
$$\Delta_1(y)=V^*(1\ot y)V.$$ 
Then $\Delta_1$ is a coproduct on the von Neumann algebra $N$ and its restriction to $Q$ is a coproduct on the C$^*$-algebra. It extends the original coproduct on $B$.
\eprop

The proof is completely similar as in Proposition \ref{prop:2.10a}. We can use Proposition \ref{prop:2.29} combined with Proposition \ref{prop:4.7a}. We could also use that $\Delta(b)=V^*(1\ot b)V$ is true in the algebra $B\ot AB$ and therefore also holds in the $^*$-representation on $\mathcal H\ot \mathcal H$.
\ssnl
Again, as for the coproduct on $M$, also here we will in the sequel use $\Delta$ for the extended coproduct on $N$.
\ssnl
This is in agreement with the definition of the coproduct on the dual we find from the general theory of locally compact quantum groups, see e.g. Definition 4.5 in \cite{VD-sigma} and remember that the coproduct on the dual in the theory of locally compact quantum groups is flipped.
%\ssnl
\keepcomment {More references, internal, external? More arguments?} %{}
\nl
\bf From the left integral from the algebra $B$ to the left Haar weight on $N$ \rm
\nl
We now have to look at the weights on the dual. For the construction of the extension of $\widehat\varphi$ from $B$ to the operator algebras, in the spirit of what we have done for $\psi$ before in this section, we need the following analogue of Proposition \ref{prop:lha}.

\prop\label{prop:4.9a}
The subspace $\widehat\Lambda(B)$ is made into a left Hilbert algebra with the $^*$-algebra structure inherited from $B$. The associated von Neumann algebra is precisely $N$.
\eprop

\bew
First observe that the map $\widehat\Lambda$ is injective with dense range in $\mathcal H$. This follows from the fact that $\Lambda:A\mapsto\mathcal H$ is injective and has dense range, together with the fact that $a\mapsto \psi(S(\,\cdot\,)a)$ is bijective from $A$ to $B$.
\ssnl
Then, as in the proof of Proposition \ref{prop:lha}, we just have to argue that the involution is a preclosed operator. To prove this we show that
\begin{equation*}
\langle\widehat\Lambda(b^*),\Lambda(x)\rangle=\langle \Lambda(S(x)^*),\widehat\Lambda(b)\rangle
\end{equation*}
for all $x\in A$ and $b\in B$. This implies the result. 
\ssnl
To prove this formula we use that, given $b=\psi(S(\,\cdot\,)a)$, 
\begin{align*}
\langle\widehat\Lambda(b^*),\Lambda(x)\rangle
&=\langle S(x)^*,b^*\rangle 
=\langle x,b\rangle^- \\
&=\psi(S(x)a)^-
=\psi(a^*S(x)^*)\\
&=\langle\Lambda(S(x)^*),\Lambda(a)\rangle
=\langle \Lambda(S(x)^*),\widehat\Lambda(b)\rangle.
\end{align*}
\ebew

Now we can proceed as in the case of the construction of the right Haar weight $\psi$ on the pair $(M,\Delta)$. We first extend the integral $\widehat\varphi$ from $B$ to a normal semi-finite weight, still denoted by $\widehat\varphi$ on the von Neumann algebra $N$. It is obtained as the weight constructed from the left Hilbert algebra $\widehat\Lambda(B)$. It is a left Haar weight.
\ssnl
\oldcomment{
We use the formula obtained in Proposition \ref{prop:JIV} in Section \ref{sect:hwt} to define the unitary antipode $\widehat R$ on $N$ by $\widehat R(y)=Jy^*J$. This will flip the coproduct on $N$. The composition $\widehat\varphi\circ \widehat R$
 will provide the right Haar weight on $N$.
 \ssnl}
 Therefore we obtain the dual locally compact quantum group.
 
 \stel
 The pair $(N,\Delta)$ where $N$ is the von Neumann algebra generated by $B$ on the Hilbert space $\mathcal H$ and where $\Delta$ is the extension of the coproduct from $B$ to $N$ as Proposition \ref{prop:2.11}, is a locally compact quantum group.
 \estel
 
 Because the Haar weigths are unique, we obtain the same weights as the ones we would obtain by applying the general procedure. 
 \ssnl
 The dual left Haar weight in the general theory is also constructed as the weight associated with the a left Hilbert algebra. To finish the comparison of the two methods, we will just show that the left Hilbert algebra from Proposition \ref{prop:4.9a} we have used here is a subalgebra of the left Hilbert algebra used in the general theory.
\ssnl
We refer to Proposition 4.10 in Section 4 of \cite{VD-sigma}. There we have the following result. We use the notations of that paper.

\prop
Denote by $\mathfrak{\mathcal N}$ the subspace of elements $y\in\widehat M$ of the form $(\omega\ot\iota)W$ where $\omega$ is a normal linear functional on $M$ such that there is a vector, denoted by $\widehat\Lambda(y)$, satisfying $\langle \widehat\Lambda(y),\Lambda(x)\rangle=\omega(x^*)$ for all $x\in \mathfrak N_\varphi$. The space  $\widehat \Lambda(\widehat{\mathcal N} \cap \widehat{\mathcal N}^*)$ is a left Hilbert algebra.
\eprop

In order to translate this result to our framework, observe that the coproduct $\Delta$ on $\widehat M$ is given by $\Delta(y)=W(y\ot 1)W*$ (see Lemma 4.4 in \cite{VD-sigma} whereas here we have $\Delta$ on $N$ given by $\Delta(y)=V^*(1\ot y)V$ (see Proposition \ref{prop:2.11}). Therefore, to translate the above proposition to our framework, we need to replace $W$ by $V^*$ and use the flip map. The result is the following, now using the notations from this paper.

\prop
Denote by $\mathfrak{\mathcal N}$ the subspace of elements $y\in N$ of the form $(\iota\ot \omega)V^*$ where $\omega$ is a normal linear functional on $M$ such that there is a vector, denoted by $\widehat\Lambda(y)$, satisfying $\langle \widehat\Lambda(y),\Lambda(x)\rangle=\omega(x^*)$ for all $x\in \mathfrak N_\varphi$. The space  $\widehat \Lambda(\widehat{\mathcal N} \cap \widehat{\mathcal N}^*)$ is a left Hilbert algebra.
\eprop

Finally we show that this left Hilbert algebra contains $\widehat\Lambda(B)$.

\prop
Let $b=\psi(S(\,\cdot\,)a)$ with $a\in A$. Denote by $\omega_1$ the functional $\psi(\,\cdot\,a)$. Then $(\iota\ot\omega)V^*=\gamma(b)$. In particular, the map $\widehat\Lambda(b)=\Lambda(a)$  where now $\widehat\Lambda(b)$ is as defined in the previous proposition.
\eprop
\bew
We use that $V^*=(\iota\ot S)V$. Then we have, for all $x\in A$, 
\begin{align*}
((\iota\ot\omega_1)V^*)\Lambda(x)
&=\sum_{(x)} \Lambda(x_{(1)})\omega_1(S(x_{(2)}))\\
&=\sum_{(x)} \Lambda(x_{(1)})\omega(x_{(2)})\\
&=\gamma(b)\Lambda(x).
\end{align*}
\ebew
  
 There is now still one argument missing. We have need  to show that right bounded elements for the left Hilbert algebra $\widehat\Lambda(B)$ we have in Proposition \ref{prop:4.9a} are the same as for the possibly bigger left Hilbert algebra defined in the framework of duality in \cite{VD-sigma}. This is the same as showing that any element $\xi$ in the bigger one, can be approximated with elements $\xi_n$ in the smaller one,  in such a way that $\pi_\ell(\xi_n)$ remains bounded. We know that this has to be true, but it may not be simple to carefully argue directly. We will not do that here.
 
\section{\hspace{-17pt}. Conclusions and further remarks} \label{sect:concl} % \input artikel3.tex %\newpage

In the the first paper on the subject, entitled \emph{ Algebraic quantum groups and duality I} \cite{VD-part1} we have worked with an algebraic quantum group and its dual.  In that paper, we use the term for a regular multiplier Hopf algebra. Most of the results we have obtained are available, but scattered in the literature. Therefore, this paper is more of a expository nature.
\ssnl
In this paper we treated $^*$-algebraic quantum groups.  These are now multiplier Hopf $^*$-algebras with positive integrals. 
\ssnl
In the first place we have considered some of the results from \cite{VD-part1} and added the extra information. The first key result is that a positive left integral exists if and only if a positive right integral exists. The scaling constant is always $1$. We also have an analytic structure. Among other things, we could define the square root of the modular element $\delta$.
\ssnl
The duality behaves nicely. When there are positive integrals on $A$, we also have positive integrals on the dual $\widehat A$. They relate to each other as in Plancerel's formula for locally compact abelian groups. 
\ssnl
In the next three sections of the paper, we describe in detail how an $^*$-algebraic quantum group gives rise to a locally compact quantum group in the sense of Kustermans and Vaes. One has to make use of elements of the theory of left Hilbert algebras. We have included them in two appendices.
\ssnl
Also the material treated in this paper is mostly found at other places in the literature, but in this case, we do have a different and simpler treatment than found in the original papers dealing with this problem.
\ssnl
One of the (minor) differences lies in the fact that we take the right integral as a starting point, thus working with the right regular representation. In most treatments, the left integral is used, but there are some advantages when starting with the right integral. 
 \ssnl
This work suggest a further study. We have discussed the analytic structure very briefly in Section \ref{sect:prel} of this paper. It is instructive to see how these properties appear as results for the associated locally compact quantum group as we obtain it from the algebraic quantum group in the Sections \ref{sect:lcp}, \ref{sect:hwt} and \ref{sect:dual}. We plan to do this in a forthcoming paper \cite{VD-part3}, entitled \emph{Algebraic quantum groups and duality III. The analytic structure on the operator algebras.}
\ssnl
Finally, this is all part of a bigger project, giving a comprehensive and updated treatment of the theory of locally compact quantum groups. This is planned in \cite{VD-notes}.

%%%%%%%%%% Appendices  %%%%%%%%%%%%%%%%%%%%%%%%%%

% Zaken aan te passen voor het geval we met appendices werken.

% Het volgende verandert de nummering van de sections naar hoofdletters voor de appendices

\renewcommand{\thesection}{\Alph{section}} 

\setcounter{section}{0}

% Het volgende dient voor een betere spatiëring als gevolg van een bredere A.1. bijvoorbeeld.

% \input heading-appendix.tex % Nieuwe instructies voor de appendices

\renewenvironment{stelling}{\begin{itemize}\item[ ]\hspace{-28pt}\bf Theorem \rm }{\end{itemize}}
\renewenvironment{propositie}{\begin{itemize}\item[ ]\hspace{-28pt}\bf Proposition \rm }{\end{itemize}}
\renewenvironment{lemma}{\begin{itemize}\item[ ]\hspace{-28pt}\bf Lemma \rm }{\end{itemize}}

\section{\hspace{-17pt}. Appendix. Weights and left Hilbert algebras}\label{sect:appA}  % \input artikel8a.tex %\newpage

In this appendix, we describe how a normal, faithful semi-finite weight is obtained from a left Hilbert algebra. This is used in Section  \ref{sect:lcp} to construct the right Haar weight for the locally compact quantum group obtained as the completion of an algebraic quantum group.
\ssnl
First we recall the definition of a left Hilbert algebra. We refer to  Chapter VI on left and right Hilbert algebras in \cite{T3}. The definition below is Definition 1.1 in Chapter VI of \cite{T3}.

\defin\label{defin:lha}
Let $\mathfrak A$ be a $^*$-algebra with the involution denoted as $\xi\mapsto\xi^\sharp$. Assume that there is a scalar product defined on $\mathfrak A$. Denote by $\mathcal H$ the Hilbert space completion of $\mathfrak A$ for this scalar product and consider $\mathfrak A$ as a subspace of $\mathcal H$. 
Then $\mathfrak A$ is called a {\it left Hilbert algebra} if the following conditions are satisfied:
\begin{itemize}[noitemsep]
\item[i)]
For every element $\xi\in \mathfrak A$ there is a bounded linear operator $\pi_\ell(\xi)$ on $\mathcal H$ satisfying $\pi_\ell(\xi)\eta=\xi\eta$ for all $\eta$ in $\mathfrak A$.
\item[ii)] 
For all $\xi, \eta, \zeta\in \mathfrak A$ we have $\langle \xi\eta,\zeta\rangle=\langle \eta,\xi^\sharp\zeta\rangle$.
\item[iii)] The subspace $\mathfrak A^2$ is still a dense subspace of the Hilbert space $\mathcal H$.
\item[iv)] The map $\xi\mapsto \xi^\sharp$ is preclosed as a conjugate linear operator on $\mathcal H$.
\end{itemize}
\edefin

The closure of the $^\sharp$-map is denoted with the same symbol and $\eta\mapsto \eta^\flat$ is used to denote the adjoint of $\xi\mapsto \xi^\sharp$. The domains of these maps are denoted by $\mathcal D^\sharp$ and $\mathcal D^\flat$. They are dense subspaces of the Hilbert space $\mathcal H$.

\opm
In can be shown  that the $^\sharp$-map is also the closure of its restriction to $\mathfrak A^2$. This is proven in Lemma 1.15 of Chapter VI  of [T]. It requires some work to prove this. In our approach to left Hilbert algebras \cite{VD-new}, we took the closure of the sharp map on the square $\mathfrak A^2$ from the very beginning. Furthermore, for the left Hilbert algebra we use  in Section \ref{sect:lcp} we have $\mathfrak A^2=\mathfrak A$. 
\eopm
Therefore, the following extra assumption is not a real restriction for our purpose here while it makes the treatment in this appendix somewhat simpler.

\voorw
In what follows, we assume the density of $\mathfrak A^2$ in $\mathfrak A$ for the $^\sharp$-norm. In other words, if we have vectors $\eta$ and $\eta'\in \mathcal H$ so that $\langle \xi^\sharp\xi',\eta\rangle=\langle \eta',{\xi'}^\sharp\xi\rangle$ for all $\xi,\xi'\in\mathfrak A$, then $\eta\in\mathcal D^\flat$ and $\eta^\flat=\eta'$.
\evoorw
Remember that we are working with conjugate linear operators. For the adjoint $^\flat$ of $^\sharp$ we have the basic equality
\begin{equation*}
\langle\xi_1^\sharp,\xi_2\rangle=\langle\xi_1,\xi_2^\flat\rangle^-=\langle \xi_2^\flat,\xi_1\rangle
\end{equation*}
when $\xi_1\in\mathcal D^\sharp$ and $\xi_2\in\mathcal D^\flat$. %We have made this comment in the section before already.
\nl 
From the conditions i), ii) and  iii) in Definition \ref{defin:lha} we know that the $\sigma$-weak closure of the elements $\pi_\ell(\xi)$ with $\xi\in\mathfrak A$ is  a von Neumann algebra.

\defin
We use $M$ for the von Neumann algebra generated by the operators $\pi_\ell(\xi)$ with $\xi\in\mathfrak A$. 
\edefin

\nl
\bf Right bounded elements and $\mathfrak A'$ \rm
\nl

For the construction of the associated weight on $M$ we will need right bounded vectors, see Definition 1.7 in Chapter VI of\cite{T3}.

\defin\label{defin:rbe}
A vector $\eta\in \mathcal H$ is called {\it right bounded} if there is a bounded operator on $\mathcal H$, necessarily unique and denoted by $\pi_r(\eta)$, given by $\pi_r(\eta)\xi=\pi_\ell(\xi)\eta$ for all $\xi\in \mathfrak A$. The uniqueness follows from the density of $\mathfrak A$ in $\mathcal H$. The set of right bounded elements is denoted by $\mathcal B$. 
\edefin

This definition is all we need to define the associated weight on the von Neumann algebra.

\defin\label{defin:A.4}
Define $\psi(x)$ for $x\in M$ and $x\geq 0$ by
\begin{equation}
\psi(x)=\sup\{\langle x\eta,\eta\rangle \mid \eta \in \mathcal B \text{ and }\|\pi_r(\eta)\|\leq 1\}.\label{eqn:A.1}
\end{equation}
\edefin

\oldcomment{We can just use $\psi$ in this appendix.
\ssnl}{}

We will show later that $\psi$ is a faithful normal semi-finite weight, obtain further properties of this weight and prove that $\psi(x^*x)=\langle \xi,\xi\rangle$ when $\xi\in\mathfrak A$ and $x=\pi_\ell(\xi)$ (see Theorem \ref{stel:main} below).  In order to do this, we need properties of the set of right bounded elements.

\prop\label{prop:A.7}
Let $\eta\in\mathcal B$. The operator $\pi_r(\eta)$ belongs to the commutant $M'$ of $M$. Moreover if $y\in M'$, we have $y\eta\in \mathcal B$ and $\pi_r(y\eta)=y\pi_r(\eta)$. In particular, the set $\pi_r(\mathcal B)$ is a left ideal of $M'$.
\eprop

\bew
 \oldcomment{Is het wel nodig dit bewijs te geven? \ssnl}{}
i) Let $\xi\in\mathfrak A$ and $\eta\in \mathcal B$. For all $\xi_1\in\mathfrak A$ we have
\begin{equation*}
\pi_r(\eta)\pi_\ell(\xi)\xi_1
=\pi_r(\eta)(\xi\xi_1)
=\pi_\ell(\xi\xi_1)\eta
=\pi_\ell(\xi)\pi_\ell(\xi_1)\eta
=\pi_\ell(\xi)\pi_r(\eta)\xi_1
\end{equation*}
By the density of $\mathfrak A$ in $\mathcal H$ we have $\pi_r(\eta)\pi_\ell(\xi)=\pi_\ell(\xi)\pi_r(\eta)$. Then  $\pi_r(\eta)\in M'$ because the operators $\pi_\ell(\xi)$ with $\xi\in\mathfrak A$ are dense in $M$. 
\ssnl
ii) Now let again $\eta\in\mathcal B$ and $y\in M'$. For all $\xi\in \mathfrak A$ we have
\begin{equation*}
\pi_\ell(\xi)y\eta=y\pi_\ell(\xi)\eta=y\pi_r(\eta)\xi
\end{equation*}
and it follows that $y\eta\in \mathcal B$ and $\pi_r(y\eta)=y\pi_r(\eta)$.
\ebew

To obtain more properties we need to consider the following subspace of $\mathcal B$.

\defin
We denote by $\mathfrak A'$ the space of elements $\eta\in \mathcal B$ that also belong to the domain $\mathcal D^\flat$.
\edefin

Again using the basic axioms, we find the following.

\prop
When $\eta\in\mathfrak A'$, then also $\eta^\flat\in\mathfrak A'$ and $\pi_r(\eta^\flat)=\pi_r(\eta)^*$.
\eprop
\bew
Let $\xi\in \mathfrak A$ and $\eta\in\mathfrak A'$. For all $\xi_1\in\mathfrak A$ we have
\begin{align*}
\langle \xi_1,\pi_\ell(\xi)\eta^\flat \rangle
&=\langle \pi_\ell(\xi)^*\xi_1,\eta^\flat\rangle
=\langle \xi^\sharp\xi_1,\eta^\flat\rangle\\
&=\langle \eta, \xi_1^\sharp\xi\rangle
=\langle \eta, \pi_\ell(\xi_1)^*\xi\rangle\\
&=\langle\pi_\ell(\xi_1)\eta,\xi\rangle
=\langle\pi_r(\eta)\xi_1,\xi\rangle\\
&=\langle \xi_1,\pi_r(\eta)^*\xi\rangle.
\end{align*}
It follows that $\pi_\ell(\xi)\eta^\flat=\pi_r(\eta)^*\xi$ for all $\xi$. Therefore $\eta^\flat$ is again right bounded and $\pi_r(\eta^\flat)=\pi_r(\eta)^*$.
Because it belongs to $\mathcal D^\flat$ we have $\eta^\flat\in\mathfrak A'$. 
\ebew

In fact, because we take for the $^\sharp$-map the closure of its restriction to $\mathfrak A^2$, we can easily obtain some kind of converse at this level.

\prop
Suppose that  $\eta, \eta'$ are right bounded vectors and $\pi_r(\eta)^*=\pi_r(\eta')$, then $\eta\in\mathcal D^\flat$ and $\eta^\flat=\eta'$.
\eprop
\bew
Take $\xi,\xi'\in\mathfrak A$. Then
\begin{align*}
\langle \xi^\sharp \xi',\eta \rangle
&=\langle \xi',\pi_\ell(\xi)\eta \rangle
=\langle \xi',\pi_r(\eta)\xi\rangle\\
&=\langle \pi_r(\eta')\xi',\xi\rangle
=\langle \pi_\ell(\xi')\eta',\xi\rangle\\
&=\langle \eta',(\xi')^\sharp\xi\rangle
\end{align*}
and we have $\eta\in\mathcal D^\flat$ and $\eta^\flat=\eta'$.
\ebew

\prop
The space $\pi_r(\mathfrak A')$ is a $^*$-subalgebra of $M'$. 
\eprop
\bew
Take $\xi\in \mathfrak A$ and $\eta,\eta_1\in \mathfrak A'$. Then
\begin{align*}
\langle  \xi,\pi_r(\eta)\eta_1\rangle
&=\langle \pi_r(\eta)^* \xi,\eta_1\rangle
=\langle \pi_r(\eta^\flat)\xi,\eta_1\rangle\\
&=\langle \pi_\ell(\xi)\eta^\flat,\eta_1\rangle
=\langle \eta^\flat,\pi_\ell(\xi)^*\eta_1\rangle\\
&=\langle \eta^\flat,\pi_\ell(\xi^\sharp)\eta_1\rangle
=\langle \eta^\flat,\pi_r(\eta_1)\xi^\sharp)\rangle\\
&=\langle \pi_r(\eta_1)^*\eta^\flat,\xi^\sharp)\rangle
=\langle \pi_r(\eta_1^\flat)\eta^\flat,\xi^\sharp\rangle.
\end{align*}
It follows that $\pi_r(\eta)\eta_1$ belongs to $\mathcal D^\flat$ and that $(\pi_r(\eta)\eta_1)^\flat=\pi_r(\eta_1^\flat)\eta^\flat$. Hence $\pi_r(\eta_1^\flat)\eta^\flat\in\mathfrak A'$ and because 
\begin{equation*}
\pi_r(\pi_r(\eta_1^\flat)\eta^\flat)=\pi_r(\eta_1^\flat)\pi_r(\eta^\flat)
\end{equation*}
we get that $\pi_r(\mathfrak A')$ is a $^*$-subalgebra of $M'$.
\ebew
\oldcomment{We may need to review this part another time.}

\nl
\bf There are enough right bounded elements \rm
\nl
We now continue with proving that there are \emph{enough} right bounded elements.

\prop
Let $\eta\in\mathcal D^\flat$. Define an operator $a$ on $\mathfrak A$ by $a\xi=\pi_\ell(\xi) \eta$. This operator has an adjoint, defined on $\mathfrak A$ by  $\xi \mapsto \pi_\ell (\xi)\eta^\flat$.
\eprop

\bew
For all $\xi, \xi'$ in $\mathfrak A$ we have
\begin{align*}
\langle a\xi,\xi'\rangle
&=\langle \pi_\ell(\xi)\eta,\xi'\rangle
=\langle \eta, \xi^\sharp\xi'\rangle\\
&=\langle {\xi'}^\sharp\xi,\eta^\flat\rangle
=\langle \xi,\pi_\ell(\xi')\eta^\flat\rangle.
\end{align*}
It follows that $\xi'$ belongs to the domain of the adjoint $a^*$ of $a$ and that $a^*(\xi')=\pi_\ell(\xi')\eta^\flat$.
\ebew

We use $\pi_r(\eta)$ for the closure of this operator $a$. 

\prop
For all $\eta\in\mathcal D^\flat$ we have that $\pi_r(\eta)$ is affiliated with $M'$.
\eprop

\bew
i) Define $a\xi=\pi_\ell(\xi)\eta$ on $\mathfrak A$. Then, for all $\xi,\xi_1\in\mathfrak A$ we have
\begin{equation*}
a\pi_\ell(\xi_1)\xi=a(\xi_1\xi)=\pi_\ell(\xi_1\xi)\eta=\pi_\ell(\xi_1)\pi_\ell(\xi)\eta=\pi_\ell(\xi_1)a\xi.
\end{equation*}
\ssnl
ii) Now let $\xi$ be in the domain of $a^*$ and $\xi_1,\xi_2\in \mathfrak A$. Then
\begin{align*}
\langle a\xi_1,\pi_\ell(\xi_2)\xi\rangle
&=\langle \pi_\ell(\xi_2^\sharp)a\xi_1,\xi\rangle \\
&=\langle a\pi_\ell(\xi_2^\sharp)\xi_1,\xi\rangle \\
&=\langle \pi_\ell(\xi_2^\sharp)\xi_1,a^*\xi\rangle\\
&=\langle\xi_1, \pi_\ell(\xi_2)a^*\xi\rangle.
\end{align*}
iii) Because $\pi_\ell(\mathfrak A)$ is dense in $M$ we also get
\begin{equation*}
\langle a\xi_1,x\xi\rangle=\langle \xi_1,xa^*\xi\rangle
\end{equation*}
for all $x\in M$, $\xi_1\in\mathfrak A$ and $\xi$ in the domain of $a^*$. Then $x\xi$ is in the domain of $a^*$ and $a^*x\xi=xa^*\xi$. This holds for all $\xi$ in the domain of $a^*$ and all $x\in M$. It follows that $a^*$ is affiliated with $M'$. The also $\pi_r(\eta)$ is affiliated with $M'$ because that is $a^{**}$.

\ebew

\keepcomment{This seems to be  a better and more logical argument than the one given in Lemma 1.10 and 1.11 in Chapter VI of [T]. But we have to check it carefully and refer to the literature.}{}

\prop\label{prop:1.12}
The space of right bounded vectors is dense in $\mathcal H$.  
\eprop

\bew
i) Consider $\eta$ in $\mathcal D^\flat$, the operator $\pi_r(\eta)$ and its polar decomposition $hu$. Let $f$ be a continuous function with compact support in the set of positive real numbers (including 0) and consider $f(h)$. Then $f(h)$ is well-defined in $M'$ because $\pi_r(\eta)$ is affiliated with $M'$. For all $\xi\in \mathfrak A$ we have
\begin{equation*}
\pi_\ell(\xi)f(h)\eta=f(h)\pi_\ell(\xi)\eta=f(h)\pi_r(\eta)\xi.
\end{equation*}
Because $f$ has compact support and $f(h)\pi_r(\eta)=f(h)hu$, we get a bounded operator and so $f(h)\eta$ is a right bounded vector.
\ssnl
ii) Now take for $f_n$ a function with range in $[0,1]$ and $f_n(t)=1$ for $t\in[0,n]$. Then $f_n(h)\xi\to \xi$ for all vectors $\xi$ (as this is true for a dense set of vectors in the domain of $\pi_r(\eta)$. It follows that $f_n(h)\eta\to\eta$ and this proves the density of the set of right bounded vectors.
\ebew

We now want to show that also $\mathfrak A'$ is still large enough.

\prop 
We have that $\mathfrak A'$ dense in $\mathcal H$ and that $\pi_r(\mathfrak A')$ is a non-degenerate $^*$-subalgebra of $M'$.
\eprop 

\bew
i) 
We first show that $\mathfrak A'$ is dense in $\mathcal H$. We go back to the proof of Proposition \ref{prop:1.12}.
%\ssnl
We have 
\begin{equation*}
\pi_r(f(h)\eta)^*=\pi_r(\eta)^*f(h)^*=\pi_r(\eta^\flat)f(h)=f(k)\pi_r(\eta^\flat)=\pi_r(f(k)\eta^\flat)
\end{equation*}
where $k=(\pi_r(\eta)^*\pi_r(\eta))^\frac12$. It follows that $f(h)\eta$ belongs to  $\mathcal D^\flat$ and that $(f(h)\eta)^\flat=f(k)\eta^\flat$. Then we can conclude that $f_n(h)\eta\to \eta$ and $(f_n(h)\eta)^\flat \to \eta^\flat$ and so $\mathfrak A'$ will be dense in $\mathcal H$.
\ssnl
ii) We have $\pi_r(\eta)\xi=\pi_\ell(\xi)\eta$. By the density of $\mathfrak A'$ and the non-degeneracy of $\pi_\ell(\mathfrak A)$ we have that the space $\pi_r(\mathfrak A')\mathcal H$ is dense. This proves the result.
\ebew

\keepcomment{To be completely sure, we have to check this again more carefully.}{}
\nl
\bf Left bounded elements \rm
\nl
We also will need the notion of left bounded elements. 

\defin\label{defin:A.16a}
A vector $\xi\in\mathcal H$ is called \emph{left bounded} if there is a bounded linear map $x$ on $\mathcal H$ satisfying   $x\eta=\pi_r(\eta)\xi$ for all right bounded vectors $\eta$. If it exists, it is unique because the set of right bounded vectors is dense. When $\xi\in\mathfrak A$, we have that this map is given by $\pi_\ell(\xi)$ and so $\xi$ is left bounded. Moreover, we can still use $\pi_\ell(\xi)$ to denote this linear map for every left bounded element. 
\edefin

Properties of left bounded vectors are very similar to the ones of right bounded vectors. We collect the most important ones here, without a proof.

\prop
i) For every left bounded vector $\xi$, we have $\pi_\ell(\xi)\in M$.\\
ii) For any $x\in M$ and $\xi$ left bounded, we have $x\xi$ left bounded and $\pi_\ell(x\xi)=x\pi_\ell(\xi)$.
\eprop

We see that the set of operators $\pi_\ell(\xi)$ with $\xi$ left bounded is a left ideal of $M$.

\notat\label{notat:A.18}
We use $\mathfrak A''$ for the set of left bounded vectors $\xi$ that belong to $\mathcal D^\sharp$.
\enotat

It is clear that $\mathfrak A\subseteq \mathfrak A''$ because $\mathfrak A\subseteq \mathcal D^\sharp$. Then we have the following properties.

\prop
If $\xi\in \mathfrak A''$ then $\xi^\sharp\in \mathfrak A''$ and $\pi_\ell( \xi^\sharp)=\pi_\ell(\xi)^*$. The space $\pi_\ell(\mathfrak A'')$ is a $^*$-subalgebra of $M$. 
\eprop

Because $\mathfrak A\subseteq \mathfrak A''$, we have enough elements in $\mathfrak A''$ and in particular, we have enough left bounded elements. We do not need the arguments that were necessary to get this  result for right bounded elements. We will later also show that $\mathfrak A$ is dense in $\mathfrak A''$ in the sense that we can approximate any element $\xi$ in $\mathfrak A$ by elements $\zeta$ in $\mathfrak A$ satisfying $\|\pi_\ell(\zeta)\|\leq \|\pi_\ell(\xi)\|$. Similarly for left bounded elements. See Proposition \ref{prop:A.23} further in this appendix and Theorem \ref{stel:B.4} in Appendix \ref{sect:appB}.

\oldcomment{\rood Do we need to formulate density properties here? Or can this wait?}

\nl
\bf The associated weight \rm
\nl
We  now obtain the main result of this appendix.

\stel\label{stel:main}
The map $\psi$ defined on elements $x\in M_+$ by 
\begin{equation*}
\psi(x)=\sup\{\langle x\eta,\eta\rangle \mid \eta \in \mathcal B \text{ and }\|\pi_r(\eta)\|\leq 1\}
\end{equation*}
(as in Definition \ref{defin:A.4}) 
is a faithful normal semi-finite weight on the von Neumann algebra $M$. We also have $\psi(x^*x)=\langle \xi,\xi\rangle$ when $\xi\in\mathfrak A$ and $x=\pi_\ell(\xi)$.
\estel

\bew
i) It is clear that $\psi(\lambda x)=\lambda\psi(x)$ for all $x\in M_+$ and $\lambda>0$. It is also clear that $\psi(0)=0$. In Proposition \ref{prop:B.3} of Appendix \ref{sect:appB} we  show that $\psi$ is additive.
\ssnl
ii) The weight is normal as it is the supremum of normal functionals. 
\ssnl
iii) The weight is faithful because, if $\langle x\eta,\eta\rangle=0$ for all $\eta$ we must have $x=0$ as the set of right bounded vectors is dense in $\mathcal H$.
\ssnl
iv) If $\xi\in \mathfrak A$ and $x=\pi_\ell(\xi)$ we have $\psi(x^*x)=\langle \xi,\xi\rangle$ because 
\begin{equation*}
\sup \{\pi_r(\eta)^*\pi_r(\eta) \mid  \eta \text{ is right bounded and }\|\pi_r(\eta)\|\leq 1\}=1.
\end{equation*}
This follows from the fact that $\pi_r(\mathcal A')$ is a non-degenerate $^*$-subalgebra of bounded operators on $\mathcal H$.
\ssnl
v) Consequently, $\psi$ is semi-finite.
\ebew
Next we want to show that the GNS space of the weight is the same as the original Hilbert space. We use an argument that we find in our lecture notes.

\prop\label{prop:A.21}
We have a linear map $\Lambda:\mathfrak N_\psi\to \mathcal H$ such that $\psi(x^*x)=\langle \Lambda(x),\Lambda(x)\rangle$ and $\pi_\psi(x)\eta=\pi_r(\eta)\Lambda(x)$ for all right bounded vectors $\eta$.
\eprop
\bew
i) Denote by $\mathcal K$ the subspace of $\mathcal H$ spanned by the vectors $\pi_r(\eta)^*\eta_1$ where $\eta_1\in\mathcal H$ and $\eta\in \mathcal B$. Suppose that $x\in M$ and define a linear functional on $V$ by $\rho(\pi_r(\eta)^*\eta_1)=\langle \eta_1,x\eta\rangle$. We claim that it is well-defined. To show this assume we have a sum $\sum \pi_r(\eta_i)^*\eta_{1i}=0$. Then 
\begin{align*}
\sum\langle\eta_ {1i}, \pi_\ell(\xi)\eta_i\rangle
&=\sum\langle \eta_ {1i},\pi_\ell(\xi)\eta_i\rangle\\
&=\sum\langle \eta_ {1i},\pi_r(\eta_i)\xi\rangle\\
&=\sum\langle \pi_r(\eta_i)^*\eta_ {1i},\xi\rangle=0.
\end{align*}
Then also $\sum\langle \eta_ {1i},x\eta_i\rangle=0$ for all $x\in M$. This proves the claim.
\ssnl
ii) Next we show that the map $\rho$ is bounded when $x\in \mathfrak N_\psi$. Indeed we have
\begin{align*}
|\rho(\pi_r(\eta)^*\eta_1)|^2
&=|\langle \eta_1,x\eta\rangle|^2\\
&\leq \|\eta\|^2\|x\eta\|^2 \langle x^*x\eta,\eta\\
&\leq \|\eta\|^2\|x\eta\|^2 \psi(x^*x)
\end{align*}
if $\|\pi_r(\eta)\|\leq 1$. In the limit $\pi_r(\eta)\to 1$ we get $|\rho(\eta_1)|^2\leq \psi(x^*x)\|\eta\|^2\|$. 
\ssnl
ii) Hence we have a vector $\Lambda(x)$ satisfying 
$$\langle \eta_1,x\eta\rangle=\langle \pi_r(\eta)^*\eta_1,\Lambda(x)\rangle=\langle\eta_1,\pi_r(\eta)\Lambda(x)$$ 
and therefore $x\eta =\pi_r(\eta)\Lambda(x)$. This holds for all $x\in\mathfrak N_\psi$  and $\eta\in \mathcal B$.
\ssnl
iii) We get in the end that $\psi(x^*x)=\langle \Lambda(x),\Lambda(x)\rangle$ for all $x\in\mathfrak N_\psi$.
\ebew

As a consequence of this, we get the following result.

\prop\label{prop:A.23a}
The subset $\Lambda(\mathfrak N_\psi)$ of $\mathcal H$ is equal to the set of left bounded vectors.
\eprop
\bew
We have seen in the previous proposition that $x\eta=\pi_r(\eta)\Lambda(x)$ for all $x\in \mathfrak N_\psi$ and $\eta$ right bounded. Hence $\Lambda(x)$ is left bounded. Conversely, let $\xi$ be left bounded and $x=\pi_\ell(\xi)$. Then 
\begin{equation*}
\langle x\eta,x\eta\rangle=\langle \pi_r(\eta)\xi,\pi_r(\eta)\xi\rangle
\end{equation*}
and if $\|\pi_r(\eta)\|\leq 1$ we obtain that $\langle x\eta,x\eta\rangle\leq \|\pi_\ell(\xi)\|^2$ and so $\psi(x^*x)$ will be finite. Hence $x\in\mathfrak N_\psi$.
\ebew
\oldcomment{Is this sufficient? We have to formulate the next as a proposition. We refer to it.}{}
\nl
We will need the following. This is part of Theorem 1.26 in Chapter 1 of \cite{T3}. For the proof, we refer to Theorem \ref{stel:B.4} in Appendix \ref{sect:appB}.

\prop\label{prop:A.23}
For all $\xi$ left bounded, there is a sequence $\xi_n$ in $\mathfrak A$ so that $\xi_n\to\xi$ while $\|\pi_\ell(\xi_n)\|\leq \|\pi_\ell(\xi\|$ for all $n$.
\eprop

Here this means that, given $x\in \mathfrak N_\psi$ there is a sequence $(a_n)\in A$ such that $\Lambda(a_n)\to\Lambda(x)$ while $\|a_n\|$ remains bounded. This is used to prove the invariance of the weight in Proposition \ref{prop:2.18} in Section \ref{sect:lcp}.

\oldcomment{\rood Refer to appendix B}

\section{\hspace{-17pt}. Appendix. Miscellaneous technical results}\label{sect:appB}  % \input artikel8b.tex% \newpage

In this appendix, we treat some of the more complicated results about left Hilbert algebras. We use the notations of the previous appendix. In particular, $\mathfrak A$ is a left Hilbert algebra and $M$ is the von Neumann generated by $\pi_\ell(\mathfrak A)$.
%\nl
\keepcomment{Possible other earlier topics}{}
\nl 
\bf The additivity of the canonical weight of a left Hilbert algebra \rm
\nl
Recall the definition of the associated weight (cf.\ Definition \ref{defin:A.4}).  
\defin
Define $\psi(x)$ for $x\in M$ and $x\geq 0$ by
\begin{equation}
\psi(x)=\sup\{\langle x\eta,\eta\rangle \mid \eta \in \mathcal B \text{ and }\|\pi_r(\eta)\|\leq 1\}.
\end{equation}
\edefin
Here $\mathcal B$ is the set of right bounded vectors.
\ssnl
We  will use the following lemma to prove the additivity of $\psi$.
\oldcomment{\ssnl In the end, we may not need this lemma. We probably can use that we have a completion of an additive functional where we are applying this. \rood Think about this. \blauw We tried but it is not obvious}{}

\lem\label{lem:A.21}
The set of positive linear functionals $x\mapsto \langle x\eta,\eta\rangle$ where $\eta$ is right bounded and $\|\pi_r(\eta)\|<1$ is upwards directed. This means that given two such functionals $x\mapsto \langle x\eta_i,\eta_i\rangle$, there is another one satisfying
\begin{equation*}
\langle x\eta_i,\eta_i\rangle\leq \langle x\eta,\eta \rangle
\end{equation*}
for all $x\in M_+$.
\elem

\bew
Put $T_i=\pi_r(\eta_i)^*\pi_r(\eta_i)$. Then $0\leq T_i<1$. Define $S_i=\frac{T_i}{1-T_i}$, put $S=S_1+S_2$ and finally $T=\frac{S}{S+1}.$
We have $S_i\leq S$ and therefore 
\begin{equation*}
T_i=\frac{S_i}{S_i+1}\leq \frac{S}{S+1}=T.
\end{equation*}
We use here that 
\begin{equation*}
\frac{S_i}{S_i+1}=1-\frac{1}{S_i+1}
\tussenen
\frac{S}{S+1}=1-\frac{1}{S+1},
\end{equation*}
and that 
\begin{equation*}
\frac{1}{S_i+1}\geq \frac{1}{S+1},
\end{equation*}
because $S_i+1\leq S+1$. 
\ssnl
Next we claim that there is a right bounded vector $\eta$ so that $T=\pi_r(\eta)^*\pi_r(\eta)$. First we write 
\begin{equation*}
\eta'_i=\frac{1}{\sqrt{1-\pi_r(\eta_i)\pi_r(\eta_i)^*}}\eta_i.
\end{equation*}
Then $\eta'_i$ is right bounded because it is of the form $y\eta_i$ where $\eta_i$ is right bounded and $y\in M'$, see Proposition \ref{prop:A.7}. Moreoveer 
\begin{align*}
\pi_r(\eta'_i)^*\pi_r(\eta'_i)
&=\pi_r(\eta_i)^*\frac{1}{1-\pi_r(\eta_i)\pi_r(\eta_i)^*}\pi_r(\eta_i)\\
&=\pi_r(\eta_i)^*\pi_r(\eta_i)\frac{1}{1-\pi_r(\eta_i)^*\pi_r(\eta_i)}\\
&=\frac{T_i}{1-T_i}=S_i.
\end{align*}

Next we show  also $S$ has the form $\pi_r(\eta'')^*\pi_r(\eta'')$ for some right bounded element $\eta''$. Because 
$\pi_r(\eta'_i)^*\pi_r(\eta'_i)=S_i\leq S$ we have elements $u_i$ in $M'$ satisfying
%\begin{equation*}
$\pi_r(\eta'_i)=u_iS^\frac12$ for all $i$.
%\end{equation*}•
We have $S=S^\frac12(u_1^*u_1+u_2^*u_2)S^\frac12$ and $u_1^*u_1+u_2^*u_2$ will be the range project of $S^\frac12$.
Then we put $\eta''=u_1^*\eta'_1+u_2^*\eta'_2$. This is a right bounded vector and we have
$\pi_r(\eta'')=u_1^*u_1 +u_2^*u_2S^\frac12$. This will be equal to $S^\frac12$. 
%\ssnl
\keepcomment{We have a similar argument in the book of Takesaki, see Lemma 2.1 of Chapter VII on weight (on page 80). In fact this is a standard result about hereditary subalgebras. We could reformulate it and refer e.g. to the book of Pedersen? Lemma 5.1.2 on page162}
\ssnl
Then, as similar argument as before will give that also $T=\pi_r(\eta)^*\pi_r(\eta)$ for some right bounded elements $\eta$. \oldcomment{\rood We still have to argue this!}{}
\ssnl
ii) All together we find a right bounded element $\eta$ with $\|\pi_r(\eta)\|< 1$ and
$\langle x\eta_i,\eta_i \rangle\leq \langle x\eta,\eta \rangle$ for all positive $x$. 
\ebew

\keepcomment{\rood We should review this.}{}

\prop\label{prop:B.3}
 When $x,y$ are positive element in the von Neumann algebra $M$, then $\psi(x+y)=\psi(x)+\psi(y)$. 
\eprop

\bew
i) . To prove additivity, we take $x_1,x_2$ in $M_+$ and put $x=x_1+x_2$. We obviously have $\psi(x)\leq \psi(x_1)+\psi(x_2)$. 
\ssnl
To prove the other inequality, we may assume that $\psi(x_i)<\infty$ for $i=1,2$. Indeed, if e.g. $\psi(x_1)=\infty$ then also $\psi(x)=\infty$ because that is bigger.
\ssnl
 We choose $\varepsilon>0$ and right bounded vectors $\eta_1$ and $\eta_2$ with the property that
\begin{equation*}
\langle x_i\eta_i,\eta_i \rangle > \psi(x_i)-\varepsilon
\tussenen \|\pi_r(\eta_i)\|<1.
\end{equation*}
By the lemma we have a right bounded vector $\eta$ so that 
\begin{equation*}
\langle x\eta_i,\eta_i \rangle\leq \langle x\eta,\eta \rangle
\tussenen
\|\pi_r(\eta)\|\leq 1.
\end{equation*}
Then 
\begin{equation*}
\psi(x_1)+\psi(x_2)\leq \langle (x_1+x_2)\eta,\eta\rangle+2\varepsilon\leq \psi(x_1+x_2)+2\varepsilon.
\end{equation*}
It follows that $\psi(x_1)+\psi(x_2)\leq \psi(x_1+x_2)$. 
\ebew

\bf The density of $\mathfrak A$ in $\mathfrak A''$ \rm
\nl
Let $\mathfrak A$ be a left Hilbert algebra and $\mathfrak A''$ it completion as in Notation \ref{notat:A.18} of Appendix \ref{sect:appA}. Then we will prove the following theorem in this appendix.

\stel\label{stel:B.4}
i) For every $\xi\in \mathfrak A''$ there is a sequence $(\xi_n)$ in $\mathfrak A$ such that $\xi_n\to\xi$ and $\xi_n^\sharp\to \xi^\sharp$ and such that $\|\pi_\ell(\xi_n)\|\leq \|\pi_\ell(\xi)\|$ for all $n$. Hence $\pi_\ell(\xi_n)\to\pi_\ell(\xi)$ in the strong-$^*$ operator topology.
\ssnl
ii) For every left bounded element $\xi$ there is a sequence $(\xi_n)$ in $\mathfrak A$ such that $\xi_n\to\xi$  and such that $\|\pi_\ell(\xi_n)\|\leq \|\pi_\ell(\xi)\|$ for all $n$. Hence $\pi_\ell(\xi_n)\to\pi_\ell(\xi)$ in the strong operator topology.
\estel
 
The result is found in \cite{T3} as Theorem 1.26 in the chapter on Left Hilbert algebras and Right Hilbert algebras. We prove it in different steps. We mostly follow the different steps in the proof given in \cite{T3}. For the easier steps, we only give indications. For some steps, we have a simpler and more direct argument.
\nl
Fix an  element $\xi\in \mathfrak A''$. Choose a sequence $(\zeta_n)$ in $\mathfrak A$ such that $\zeta_n\to \xi$ and $\zeta_n^\sharp\to \xi^\sharp$. This is possible because $\xi\in \mathcal D^\sharp$ and the $^\sharp$-operator is defined as the closure of its restriction to $\mathfrak A$.
\ssnl
Denote $x=\pi_\ell(\xi)$ and $x_n=\pi_\ell(\zeta_n)$. 

\lem We have $x_n\eta\to x\eta$ and $x_n^*\eta\to x_n\eta$ for all right bounded vectors $\eta$.
\elem
\bew
For any right bounded vector $\eta$ we have
\begin{equation*}
x_n\eta=\pi_\ell(\zeta_n)\eta=\pi_r(\eta)\zeta_n\to \pi_r(\eta)\xi=\pi_\ell(\xi)\eta=x\eta
\end{equation*}
and similarly $x_n^*\eta\to x^*\eta$ because also $\zeta_n^\sharp\to\xi^\sharp$.
\ebew

The next step is a 2x2-matrix trick to get self-adjoint operators as we have in Takesaki's proof. We use $\mathcal H_1$ for the Hilbert space direct sum $\mathcal H\oplus\mathcal H$ of $\mathcal H$ with itself.

\lem
We set 
\begin{equation*}
X_n=\left(\begin{matrix} 0 & x_n^*\\ x^n & 0 \end{matrix}\right)
\tussenen
X=\left(\begin{matrix} 0 & x^*\\ x & 0 \end{matrix}\right).
\end{equation*}
Then we have self-ajoint operators on $\mathcal H_1$ such that $X_n\eta\to X\eta$ for a dense subspace $\mathcal B_1$ of vectors $\eta\in \mathcal H_1$.
\elem

The proof is obvious. For the subspace $\mathcal B_1$ we can take elements of the form $(\eta_1,\eta_2)$ where $\eta_1,\eta_2$ are right bounded in $\mathcal H$.
\ssnl
Now we take a slightly different path, but still in the same spirit of the proof of Takesaki in \cite{T3}.
\ssnl
We have the following  convergence property.

\lem
We have that
\begin{equation*}
\frac{1}{1+X_n^2}\to \frac{1}{1+X^2}
\end{equation*}
in the strong operator topology.
\elem

\bew
i) We have
\begin{align*}
\frac{1}{1+iX_n} - \frac{1}{1+iX}
&=\frac{1}{1+iX_n}((1+iX)-(1+iX_n)) \frac{1}{1+iX}\\
&=i\frac{1}{1+iX_n}(X-X_n) \frac{1}{1+iX}.
\end{align*}
Now we observe that $\|\frac{1}{1+iX_n}\|\leq 1$. If now $\eta$ is a vector in $\mathcal H_1$ such that $\frac{1}{1+iaX}\eta$ belongs to $\mathcal B_1$, we have that 
\begin{equation*}
\frac{1}{1+iX_n}\eta \to \frac{1}{1+iX}\eta.
\end{equation*}
Such vectors are still dense and because all these operators are bounded in norm by $1$, we get that 
\begin{equation*}
\frac{1}{1+iX_n}\eta \to \frac{1}{1+iX}\eta.
\end{equation*} for all vectors $\eta\in\mathcal H_1$. 
\ssnl
ii) Similarly we have that $\frac{1}{1+iX_n} \to \frac{1}{1+iX}$ in the strong operator topology. 
\ssnl
iii) Finally we have
\begin{equation*}
\frac{1}{1-iX}+\frac{1}{1+iX}=\frac{2}{1+X^2}
\end{equation*}
and the same for $X_n$. We can conclude that 
\begin{equation*}
\frac{1}{1+X_n^2}\to \frac{1}{1+X^2}
\end{equation*}
in the strong operator topology.
\ebew

Next we use the following result we find as Proposition 2.3.2 in \cite{Pe}.

\prop\label{prop:B.5}
Each continuous function $g$ on $\mathbb R$ such that  $g(0) = 0$ and $|g(t)|\leq \alpha|t|+\beta$ 
for some positive $\alpha$ and $\beta$ is strongly continuous.
\eprop

A function $g$ is called strongly continuous if for each net $(x_i)$ of self-adjoint bounded operators, with limit $x$, the net $(g(x_i))$ converges strongly to $g(x)$. We will apply this with the function $g$ defined on the interval $[0,1]$ by
\begin{equation*}
g(s)=\left\{\begin{matrix} 2s & \text{if } &0\leq s<\frac12 \\ 1 &\text{if } &\frac12\leq s\leq1\end{matrix}\right.
\end{equation*}

We get the following result.

\lem
Define $g$ as above and $f(t)=g(\frac{1}{1+t})$ for $0\leq t$. Then 
\begin{equation*}
f(x_n^*x_n)\to f(x^*x)=1
\tussenen
f(x_nx_n^*)\to f(xx^*)=1
\end{equation*}
strongly.
\elem
\bew
i) If we apply Proposition \ref{prop:B.5} we find that $f(X_n^2)\to f(X^2)$ strongly. Because  
\begin{equation*}
X^2=\left(\begin{matrix} xx^* & 0 \\ 0 &x^*x \end{matrix}\right)
\end{equation*}
and similarly for $X_n^2$ we get
\begin{equation*}
f(x_n^*x_n)\to f(x^*x)
\tussenen
f(x_nx_n^*)\to f(xx^*).
\end{equation*}
And because $\|x\|\leq 1$  we have that $\frac12\leq \frac{1}{1+X^2}\leq 1$ so that $f(X)=1$.
\ebew

For the further steps, we again follow the argument of Takesaki.
\ssnl
Now we define $\xi_n=f(x_nx_n^*)\zeta_n$. 
\lem
The elements $\xi_n$ belongs to $\mathfrak A''$ and $\xi_n^\sharp=f(x_n^\sharp x_n)\zeta_n^\sharp$. Furthermore
$\xi_n\to \xi$  while now $\|\pi_\ell(\xi_n)\|\leq 1$ for all $n$.
\elem
\bew
i) We get 
\begin{align*}
\|\xi_n-\xi\|
&\leq \|f(x_nx_n^*)(\zeta_n-\xi_n)\| + \| f(x_nx_n^*)\xi-\xi \|\\
&\leq \|\zeta_n-\xi\|+\|f(x_nx_n^*)\xi-f(xx^*)\xi\|\to 0.
\end{align*}
\ssnl
ii) Similarly we get
\begin{equation*}
\|\xi_n^\sharp-\xi^\sharp\|\to 0.
\end{equation*}
\ssnl
iii)
and also, because
\begin{align*}
&f(t^2)|t|<1  \text{ if } 0\leq |t| \leq1 \\
&f(t)|t|=\frac{2|t|}{1+t^2} \leq 1 \text { if } 1\leq |t|
\end{align*}
we get  $\|\pi_\ell(\xi_n)\|=\|f(x_nx_n^*)x_n\|\leq 1$. 
\ebew

\keepcomment{\rood We have to check this carefully.
\nl}
The elements $\xi_n$ approximating $\zeta$ are no longer expected to be in $\mathfrak A$ as we want. By approximating the function $f$ by polynomials in a careful chosen way, we can finally complete the proof of item i) of Theorem \ref{stel:B.4}.

\prop
There is a sequence $(\xi'_n)$ in $\mathfrak A$ such that $\xi'_n\to\xi$ and $\zeta_n^\sharp\to \xi^\sharp$ and such that $\|\pi_\ell(\xi'_n)\|\leq \|\pi_\ell(\xi)\|$ for all $n$. Hence $\pi_\ell(\xi'_n)\to\pi_\ell(\xi)$ in the strong-$^*$ operator topology.
\eprop

\bew
We now choose for each $n$ a polynomial $p_n$ that approximates $f$. Then we put $\xi'_n=p_n(x_nx_n^*)\zeta_n$. Again $\xi'_n\in\mathfrak A$ and $(\xi'_n)^\sharp= p_n(x_n^*x_n)\zeta_n^\sharp$.  We have now
\begin{align*}
&\|\xi'_n-\xi_n\|\leq \|p_n(x_nx_n^*)-f(x_nx_n^*)\|\,\|\zeta_n\|  \\
&\|(\xi'_n)^\sharp-\xi_n^\sharp\|\leq \|p_n(x_n^*x_n)-f(x_n^* x_n)\|\,\|\zeta_n^\sharp\|
\end{align*}
We can choose $p_n$ so that 
\begin{equation*}
\|\xi'_n-\xi_n\|\leq \frac1n 
 \tussenen
 \|(\xi'_n)^\sharp-\xi_n^\sharp\|\leq \frac1n.
\end{equation*}
Furthermore we get
\begin{align*}
\|\pi_\ell(\xi'_n\|
&=\|p_n(x_nx_n^*)x_n^*\| \\
&\leq \| p_n(x_nx_n^*)-f(x_nx_n^*)\|\|x_n\|+\|f(x_nx_n^*)x_n\|
\end{align*}
We can choose $p _n$ so that also 
$$\| p_n(x_nx_n^*)-f(x_nx_n^*)\|\|x_n\|\leq \frac1n.$$ 
Then we get $\|\pi_\ell(\xi'_n\|\leq 1+\frac 1n$. Finally we just have to scale $\xi'_n$ by the factor $1+\frac1n$.
\ebew

This proves item i) of the theorem. We get item ii) as a standard consequence of this.

\prop
For every left bounded element $\xi$ there is a sequence $(\xi_n)$ in $\mathfrak A$ such that $\xi_n\to\xi$  and such that $\|\pi_\ell(\xi_n)\|\leq \|\pi_\ell(\xi)\|$ for all $n$. Hence $\pi_\ell(\xi_n)\to\pi_\ell(\xi)$ in the strong operator topology.
\eprop
\bew
Assume that $\xi$ is left bounded and that $\|\pi_\ell(\xi)\|\leq 1$. 
We can approximate $\xi$ by elements of the form $\pi_\ell(\xi_1)\xi$ where $\xi_1\in\mathfrak A$ and where $\|\pi_\ell(\xi_1)\|\leq 1$. These elements belong to $\mathfrak A''$. We can approximate them by elements $\xi_n\in \mathfrak A$ and such that $\|\pi_\ell(\xi_n)\|\leq 1$ for all $n$.
\ebew

\keepcomment{\rood Why is this more complicated in the proof of Takesaki. Also remark why we no longer have strong$^*$-operator convergence.}{}

%\nl
\keepcomment{\rood There are a lot of misprints in this part of the book of Takesaki!}{}
\nl

%%%%%%%% References %%%%%%%%%%%%%

% \input artikel9.tex % Referenties

\end{document}